\documentclass[twocolumn,secthm]{autart}

\usepackage{graphicx}
\usepackage{cite}
\usepackage{amsmath,amssymb,amsfonts,bm}
\usepackage{mathtools}
\usepackage{algorithmic}
\usepackage{graphicx}
\usepackage{textcomp}
\usepackage{dsfont}
\usepackage{color}
\usepackage{picins}
\usepackage{afterpage}
\usepackage{epstopdf}        
\usepackage{enumerate}
\usepackage{lscape}
\usepackage{multicol}
\usepackage{multirow}
\usepackage{calligra}
\usepackage{booktabs}
\usepackage{mathtools}
\usepackage{framed}  
\usepackage{empheq}
\usepackage{graphics}
\usepackage{epsfig}
\usepackage{epstopdf}
\usepackage{enumerate}
\usepackage{graphicx} 
\usepackage{comment} 
\usepackage{apptools}
\usepackage{subcaption}
\usepackage{theoremref}
\usepackage{bbm}

\graphicspath{{./},{./figures/}}

\DeclareMathOperator*{\arginf}{arg\,inf}

\newtheorem{definition}{Definition}
\newtheorem{asmp}{Assumption}
\newtheorem{example}{Example}

\def\BibTeX{{\rm B\kern-.05em{\sc i\kern-.025em b}\kern-.08em
    T\kern-.1667em\lower.7ex\hbox{E}\kern-.125emX}}

\begin{document}

\begin{frontmatter}

\vspace{-0.4cm}
\title{On Persistently Resetting Learning Integrators:\\ A Framework For Model-Free Feedback Optimization\thanksref{footnoteinfo}} 

\thanks[footnoteinfo]{This work was not presented at any conference. This work was supported in part by the grants NSF ECCS CAREER 2305756 and AFOSR YIP: FA9550-22-1-0211. Corresponding Author: Jorge I. Poveda.}

\author{Mahmoud Abdelgalil}\ead{mabdelgalil@ucsd.edu},    
\author{Jorge I. Poveda}\ead{jipoveda@ucsd.edu},               
    
\address{Department of Electrical and Computer Engineering, University of California, San Diego, La Jolla, CA, USA.}

\begin{keyword}                           
Hybrid systems, Adaptive systems, Nonlinear Control, Optimization.             
\end{keyword}                             

\begin{abstract}
We study a novel class of algorithms for solving model-free feedback optimization problems in dynamical systems. The key novelty is the introduction of \emph{persistent resetting learning integrators} (PRLI), which are integrators that are reset at the same frequency at which the plant is dithered using exploratory signals for model-free optimization. It is shown that PRLIs can serve as core mechanisms for real-time gradient estimation in online feedback-optimization tasks where only cost function measurements are available. In particular, unlike existing approaches based on approximation theory, such as averaging or finite-differences, PRLIs can produce global real-time gradient estimates of cost functions, with uniformly bounded perturbations of arbitrarily small magnitude. In this sense, PRLIs function as robust \emph{hybrid} ``Oracles" suitable for interconnection with discrete-time optimization algorithms that optimize the performance of continuous-time dynamical plants in closed-loop operation. Compared to existing methods, PRLIs yield \emph{global} stability properties for a broad class of cost functions, surpassing the local or semi-global guarantees offered by traditional approaches based on perturbation and approximation theory. The proposed framework naturally bridges physical systems, modeled as continuous-time plants where continuous exploration is essential, with digital algorithms, represented as discrete-time optimization methods. This integration of physical and digital components allows for the application of a broad range of well-established optimization techniques, each offering unique performance benefits suited to different applications—such as accelerated convergence through momentum, safety assurances via projections, etc. Closed-loop stability guarantees are established using unifying tools from hybrid systems theory. The main results are illustrated using different numerical examples.  
\end{abstract}

\end{frontmatter}

\section{Introduction}
In many engineering applications, optimizing the steady-state output map of a dynamic plant is highly desirable. Such applications include, but are not limited to, power systems \cite{ghaffari2014power,colombino2019online}, transportation and traffic control \cite{yu2021extremum}, chemical processes \cite{guay2004adaptive,wang1999optimizing}, and robotics \cite{abdelgalil2023singularly}. When sufficient information is available regarding the model of the plant, and if its input-output behavior is well-understood and remain unchanged over time, then standard numerical optimization schemes maybe interconnected with the plant. This paradigm, of interconnecting optimization algorithms with plant dynamics, has received significant attention over the past two decades and continues to be a flourishing research field till this day \cite{colombino2019online,hauswirth2024optimization,scheinker2024100}.

For most applications, however, operating conditions are continually changing, the plant is often subjected to disturbances, and precise plant models may not be available. In these circumstances, model-based optimization methods are inapplicable and may lead to { instabilities} or poor performance. By contrast, model-free online optimization methods probe the input-output behavior of the plant in real-time. Thus, under sufficient regularity conditions, model-free methods are able to continually adapt to a changing plant.

\vspace{-0.2cm}
The majority of model-free online optimization algorithms attempt to estimate the gradient of the steady-state output function from collected measurements of its value. Indeed, gradient estimation through function evaluations is the paradigmatic approach common to all gradient-based \emph{zeroth-order optimization} methodologies. Consequently, model-free methods face the fundamental exploration-exploitation dilemma, i.e., experimenting with available options to gain new information (exploration) and utilizing known information (exploitation) to optimize the cost function. 

Abstractly, a model-free (gradient-based) optimization scheme consists of two components: i) an \emph{oracle} for gradient estimation (exploration), ii) and an update rule used as an optimization routine (exploitation). The two components of the algorithm are often mixed, especially when the optimization problem needs to be solved in real time and the cost function characterizes the input-to-output map of a dynamical system. Indeed, this is precisely the situation in \emph{extremum seeking} (ES) control \cite{Leblanc1922}, one of the most popular and successful families of model-free methods for real-time feedback-based optimization of dynamical systems. 

{ ES techniques have addressed the exploration \!/\! exploitation dilemma using different approaches, including (but not limited to) standard averaging theory \cite{krstic2000stability,ariyur2003real,tan2010extremum,oliveira2016extremum,zhu2022extremum,yilmaz2024prescribed,scheinker2024100}, sampled-data techniques \cite{khong2013unified,hazeleger2022sampled}, finite-differences-based approaches \cite{khong2015extremum}, Lie-Bracket approximation and second-order averaging techniques \cite{durr2013lie,abdelgalil2022recursive}, and averaging for hybrid dynamical systems \cite{poveda2017framework,abdelgalil2025lie}. Nevertheless, the majority of ES techniques suffer from two main issues. The first issue is related to the fact that ES control relies on two distinct time-scale separation arguments; namely, averaging and singular perturbation \cite{krstic2000stability} (see also \cite{abdelgalil2023multi} for more details on the applications of singular perturbation and averaging in model-free optimization). As a result, the design of ES techniques typically necessitates sequential tuning of two distinct parameters; the frequency of oscillation of the exploratory signals and the timescale separation constant with respect to the plant dynamics. In particular, even in the case of a static cost function, the design of an ES algorithm necessitates the tuning of the frequency of the exploration signals for every compact set of initial conditions \cite{krstic2000stability,durr2013lie,zhu2022extremum}. The second issue is the fact that, even for quadratic cost functions, ES control can lead to \emph{finite-time escape} if initialized outside the appropriate set of initial conditions. From this perspective, ES is not robust to perturbations that push the system outside of the intended set of initial conditions. As a consequence of these two issues, the majority of ES algorithm in the literature are, at best, semi-globally practically asymptotically stable \cite{teel1999semi}. Although our earlier work \cite{abdelgalil2024initialization} addressed the above issues, the assumptions on the cost function adopted in \cite{abdelgalil2024initialization} are restrictive. Moreover, the results in \cite{abdelgalil2024initialization} were concerned with gradient descent only and do not apply when a more sophisticated optimization scheme, e.g. projected gradient descent, is employed.}

{ To address the above shortcomings, we propose a novel model-free zeroth-order optimization framework}. A key innovation of our work is the use of \emph{persistently resetting learning integrators} (PRLI). These are integrators that reset in sync with the exploratory signals used for gradient estimation. We demonstrate that PRLIs can act as core mechanisms for real-time gradient estimation in online feedback-optimization tasks where only cost function measurements (zeroth-order information) are accessible. Specifically, we show that PRLIs can deliver \emph{global}, real-time gradient estimates of the cost function, with uniformly bounded perturbations of arbitrarily small magnitude. Consequently, PRLIs serve as robust \emph{oracles} suitable for coupling with discrete-time optimization algorithms that aim to optimize the performance of continuous-time dynamical plants in closed-loop operation. In contrast to established approaches--such as those relying on averaging theory--PRLIs can guarantee \emph{global stability} properties for a wide range of cost functions, thus surpassing the typically local or semi-global assurances that traditional perturbation and approximation theories provide. Furthermore, by naturally bridging continuous-time physical systems (plants) and discrete-time digital algorithms (optimization routines), the proposed framework supports a seamless integration with well-established model-based digital optimization schemes for which different computational properties are known. These include, for example, momentum-based methods for accelerated convergence and projection-based methods for safety guarantees.

\vspace{-0.3cm}
The proposed framework is most naturally stated and analyzed within the framework of Hybrid Dynamical Systems (HDS) \cite{goebel2012hybrid}, which is reviewed in Section \ref{sec:preliminaries}. In Section \ref{sec:PRLIs}, we introduce PRLIs and prove their key properties, as well as some preliminary results concerning the stability properties of the interconnection between PRLIs and a generic gradient-based optimization schemes. Section \ref{sec:stability_properties} presents the global stability properties of the proposed optimization framework for several optimization schemes and under natural conditions on the cost function. Section \ref{sec:dynamic-plants} considers the case when the cost function arises as the steady-state input-output map of a dynamical system. Section \ref{sec:numerical-results} presents numerical simulation results, and finally  Section \ref{sec:conclusion} ends with the conclusions.

\vspace{-0.3cm}
\section{Notation, Definitions, and Preliminaries}\label{sec:preliminaries}
\subsection{Notation}
The bi-linear form $\langle \cdot,\cdot\rangle:\mathbb{R}^n\times\mathbb{R}^n\rightarrow\mathbb{R}$ is the canonical inner product on $\mathbb{R}^n$, i.e. $\langle x, y\rangle = x^\top y$, for any two vectors $x,y\in\mathbb{R}^n$. We use $\|x\|=\sqrt{\langle x,x\rangle}$ to denote the $2$-norm of a vector $x\in\mathbb{R}^n$. We denote by $\mathbb{R}_{+}$ the set of positive real numbers, and by $\mathbb{R}_{\geq 0}$ the set of non-negative real numbers. For a compact set $\mathcal{A}\subset\mathbb{R}^n$ and a vector $x\in\mathbb{R}^n$, we use $|x|_{\mathcal{A}}:=\min_{\tilde{x}\in\mathcal{A}}\|x-\tilde{x}\|$. For a continuous function $V:\mathbb{R}^n\rightarrow\mathbb{R}$, we use $\mathcal{L}_c(V)$ to denote the level set $\{x\in\mathbb{R}^n~|~V(x)=c\}$ for any $c\in\mathbb{R}$. A function $\beta:\mathbb{R}_{\geq0}\times\mathbb{R}_{\geq0}\to\mathbb{R}_{\geq0}$ is of class $\mathcal{KL}$ if it is nondecreasing in its first argument, nonincreasing in its second argument, $\lim_{r\to0^+}\beta(r,s)=0$ for each $s\in\mathbb{R}_{\geq0}$, and  $\lim_{s\to\infty}\beta(r,s)=0$ for each $r\in\mathbb{R}_{\geq0}$. Throughout the paper, for two (or more) vectors $u,v \in \mathbb{R}^{n}$, we write $(u,v)=[u^{\top},v^{\top}]^{\top}$. A function $f$ is said to be $\mathcal{C}^k$ if its $k$th-derivative is continuous.
\subsection{Hybrid Dynamical Systems}
The algorithms studied in this paper will integrate both continuous-time and discrete-time dynamics. As such, they will be modeled as hybrid dynamical systems (HDS) aligned with the framework of \cite{goebel2012hybrid}, and characterized by the following equations:
\begin{align}\label{eq:hds_notation}
\mathcal{H}:~~\begin{cases}
    ~~\xi\in C, & \dot{\xi}\hphantom{^+}=F(\xi)\\
    ~~\xi\in D, & \xi^+=G(\xi),
\end{cases}
\end{align}
where $F:\mathbb{R}^n\rightarrow\mathbb{R}^n$ is called the flow map, $G:\mathbb{R}^n\rightarrow\mathbb{R}^n$ is called the jump map, $C\subset\mathbb{R}^n$ is called the flow set, and $D\subset\mathbb{R}^n$ is called the jump set. We use $\mathcal{H}=(C,F,D,G)$ to denote the \emph{data} of the HDS, and we impose the following standard assumption on $\mathcal{H}$:
\begin{asmp}\label{asmp:regularity}
The sets $C$ and $D$ are closed, and the functions $F$ and $G$ are continuous.
\end{asmp}
Any HDS that satisfies Assumption \ref{asmp:regularity} is a \emph{well-posed} HDS in the sense of \cite[Definition 6.29]{goebel2012hybrid}, see, e.g., \cite[Theorem 6.30]{goebel2012hybrid}. 

Solutions to  \eqref{eq:hds_notation} are parameterized by a continuous-time index $t\in\mathbb{R}_{\geq0}$, which increases continuously during flows, and a discrete-time index $j\in\mathbb{Z}_{\geq0}$, which increases by one during jumps. Therefore, solutions to \eqref{eq:hds_notation} are defined on \emph{hybrid time domains} (HTDs). A set $E\subset\mathbb{R}_{\geq0}\times\mathbb{Z}_{\geq0}$ is called a \textsl{compact} HTD if $E=\cup_{j=0}^{J-1}([t_j,t_{j+1}],j)$ for some finite sequence of times $0=t_0\leq t_1\ldots\leq t_{J}$. The set $E$ is a HTD if for all $(T,J)\in E$, $E\cap([0,T]\times\{0,\ldots,J\})$ is a compact HTD. A hybrid arc $\xi$ is a function defined on a HTD. In particular, $\xi:\text{dom}(\xi)\to \mathbb{R}^n$ is such that $\xi(\cdot, j)$ is locally absolutely continuous for each $j$ such that the interval $I_j:=\{t:(t,j)\in \text{dom}(\xi)\}$ has a nonempty interior. A hybrid arc $\xi:\text{dom}(\xi)\to \mathbb{R}^n$ is a solution $\xi$ to the HDS \eqref{eq:hds_notation} if $\xi(0, 0)\in C\cup D$, and:
\begin{enumerate}
\item For all $j\in\mathbb{Z}_{\geq0}$ such that $I_j$ has nonempty interior: $\xi(t,j)\in C$ for all $t\in I_j$, and $\dot{\xi}(t,j)= F(\xi(t,j))$ for all $t\in I_j$.
\item For all $(t,j)\in\text{dom}(\xi)$ such that $(t,j+1)\in \text{dom}(\xi)$: $\xi(t,j)\in D$ and $\xi(t,j+1)= G(\xi(t,j))$.
\end{enumerate}
A solution $\xi$ is said to be \emph{maximal} if it cannot be further extended. A solution $\xi$ is said to be \emph{complete} if the length of its domain is infinite, i.e., for every $T>0$ there exists $(t,j)\in\mathrm{dom}(\xi)$ such that $t+j>T$.
In this manuscript, we use the following definition concerning the limiting behavior of solutions of hybrid systems.
\begin{definition}\label{defn:omega-limit-set}\cite[Definition 6.23]{goebel2012hybrid}
    The $\omega$-limit set of the solutions of an HDS from a given set $\mathcal{K}\subset C\cup D$, denoted by $\Omega(\mathcal{K})$, is the set of all points $x\in\mathbb{R}^n$ for which there exists a sequence $\{\xi_i\}_{i=1}^\infty$ of solutions of the HDS $\mathcal{H}$ with $\xi_i(0,0)\in K$ and a sequence $\{(t_i,j_i)\}_{i=1}^\infty$ of points $(t_i,j_i)\in\mathrm{dom}(\xi_i)$ such that $\lim_{i\rightarrow\infty} t_i+j_i = \infty$ and $\lim_{i\rightarrow\infty} \xi_i(t_i,j_i) = x$.
\end{definition}
\subsection{Boundedness and Stability Notions for HDS}
The following stability notions for $\mathcal{H}$ will be relevant in the study of the algorithms considered herein.
\begin{definition}
    A compact set $\mathcal{A}\subset C\cup D$ is said to be Uniformly Globally Asymptotically Stable (UGAS) for the HDS $\mathcal{H}$ if there exists a class-$\mathcal{KL}$ function $\beta$ such that every solution $\xi$ of the HDS $\mathcal{H}$ satisfies
    \begin{align*}
        |\xi(t,j)|_{\mathcal{A}}\leq \beta(|\xi(0,0)|_{\mathcal{A}},t+j),
    \end{align*}
    for all $(t,j)\in\mathrm{dom}(\xi)$. 
\end{definition}
\begin{definition}
    The HDS $\mathcal{H}$ is said to be $\kappa$-Uniformly Globally Ultimately Bounded ($\kappa$-UGUB) with respect to a compact set $\mathcal{A}\subset C\cup D$ if there exists a class-$\mathcal{KL}$ function $\beta$ and a constant $\kappa>0$ such that every solution $\xi$ of $\mathcal{H}$ with $\xi(0,0)\in C\cup D$ satisfies
    \begin{align*}
        |\xi(t,j)|_{\mathcal{A}}\leq \beta(|\xi(0,0)|_{\mathcal{A}},t+j) + \kappa,
    \end{align*}
    for all $(t,j)\in\mathrm{dom}(\xi)$. 
\end{definition}
We remark that the definition of $\kappa$-UGUB we consider here implies the notion of uniform global ultimate boundedness defined, e.g. , in \cite[Definition 4]{subbaraman2016equivalence}. However, our definition, which emphasizes the role of the compact set $\mathcal{A}$, will be convenient to state some of our results and proofs.
We will also need the following definition of uniform global recurrence of a set for the HDS $\mathcal{H}$.  Although the property of recurrence has not been previously studied in model-free feedback optimization, it will be shown to be relevant for the algorithms examined in this paper.
\begin{definition}\cite[Definition 2]{subbaraman2016equivalence}\label{defn:recurrence}
    A set $\mathcal{O}\subset \mathbb{R}^n$ is said to be Uniformly Globally Recurrent (UGR) for the HDS $\mathcal{H}$ if there is no finite escape times for $\mathcal{H}$, and for each compact set $K$ there exists $T>0$ such that for each solution $\xi$ with $\xi(0,0)\in K$ either $t+j<T$ for all $(t,j)\in\mathrm{dom}(\xi)$ or there exists $(t,j)\in\mathrm{dom}(\xi)$ such that $t+j\leq T$ and $\xi(t,j)\in\mathcal{O}$.
\end{definition}
Finally, when the data of the HDS $\mathcal{H}$ depends on a small parameter $\varepsilon$, we will need the following stability notion.
\begin{definition}
    A compact set $\mathcal{A}\subset C\cup D$ is said to be Semi-Globally Practically Asymptotically Stable (SGPAS) as $\varepsilon\rightarrow 0^+$ for the HDS $\mathcal{H}$ if there exists a class-$\mathcal{KL}$ function $\beta$ such that for every compact set $K$ and every $\kappa > 0$, there exists $\varepsilon^*>0$ such that, for all $\varepsilon\in(0,\varepsilon^*)$, every solution $\xi$ of $\mathcal{H}$ with $\xi(0,0)\in K$ satisfies
    \begin{align*}
        |\xi(t,j)|_{\mathcal{A}}\leq \beta(|\xi(0,0)|_{\mathcal{A}},t+j) + \kappa,
    \end{align*}
    for all $(t,j)\in\mathrm{dom}(\xi)$. 
\end{definition}
{ \subsection{Unbiased Persistently Exciting Functions}
The following definition will play a key role in our paper:
\begin{defn}\label{defn:UPE}
    The function $v(\cdot):[0,1]\rightarrow\mathbb{R}^n$ is said to be an Unbiased Persistently Exciting (UPE) function if $v$ is continuous, and satisfies:
    \begin{align}\label{eq:unbiased-exploration}
    \int_0^1 v(\tau)\,\mathrm{d}\tau&= 0, & \int_0^1 v(\tau)v(\tau)^\top\,\mathrm{d}\tau&=  I.
\end{align}
\end{defn}

\vspace{-0.4cm}
{  A simple example of a multivariable UPE function is given by
\begin{align}\label{eq:sinusoidal-exploration}
    v(\tau)&= \sqrt{2}(\sin(2\pi\tau),\sin(4\pi\tau),\dots,\sin(2n\pi\tau)).
\end{align} 
Note that the class of UPE functions relates to \emph{persistently exciting} (PE) functions, which are common in adaptive control and system identification \cite{sastry2011adaptive}. The following Lemma is a useful tool to generate UPE functions from PE functions.}
\begin{lem}\label{lem:UPE-from-PE}\normalfont
    Let $\rho(\cdot):[0,1]\rightarrow\mathbb{R}^n$ be a continuous function for which there exists a constant $\alpha\in\mathbb{R}_{+}$ such that
    \begin{align*}
        \int_0^1 \rho(\tau)\mathrm{d}\tau\int_0^1 \rho(\tau)^\top\mathrm{d}\tau + \alpha I \preceq \int_0^1 \rho(\tau)\rho(\tau)^\top\,\mathrm{d}\tau.
    \end{align*}
    Then, the function $v(\cdot):[0,1]\rightarrow\mathbb{R}^n$ defined by
    \begin{align*}
        v(\tau)= Q^{-\frac{1}{2}} \left(\rho(\tau)-\bar{\rho}\right),
    \end{align*}
    is a UPE function, where $\bar{\rho}$ and $Q$ are given by
    \begin{align*}
        \bar{\rho}&= \int_0^1\rho(\tau)\mathrm{d}\tau, & Q_v &= \int_0^1 \rho(\tau)\rho(\tau)^\top\mathrm{d}\tau-\bar{\rho}\bar{\rho}^\top.
    \end{align*}

    \vspace{-0.8cm}
    \phantom{.}\hfill $\square$
\end{lem}

\vspace{-0.4cm}
We remark that the class of UPE functions is a substantial generalization of commonly used exploration signals in adaptive systems, e.g. \eqref{eq:sinusoidal-exploration}. Indeed, Definition \ref{defn:UPE} does not require periodicity in $\tau$. Moreover, the flexibility in the definition allows for the use of \emph{randomly} constructed exploration signals, e.g. signals constructed via Lemma \ref{lem:UPE-from-PE} from a (fixed) sample path of a continuous random process. The details of such a construction are out-of-scope of the current manuscript. Nevertheless, in Section \ref{sec:numerical-results}, we provide an example of a UPE function that is not of the form \eqref{eq:sinusoidal-exploration}. We now give a proof of Lemma \ref{lem:UPE-from-PE}.
\begin{pf}
    Define the function $\tilde{v}$ by
    \begin{align*}
        \tilde{v}(\tau)= \rho(\tau)-\bar{\rho},
    \end{align*}
    which is, by construction, continuous since $\rho$ is continuous. Direct computation gives
    \begin{align*}
        \int_0^1 \tilde{v}(\tau)\,\mathrm{d}\tau&= \int_0^1 (\rho(\tau)-\bar{\rho})\,\mathrm{d}\tau = 0.
    \end{align*}
    Similarly, we compute that
    \begin{align*}
        \int_0^1 \tilde{v}(\tau)\tilde{v}(\tau)^\top\,\mathrm{d}\tau &= \int_0^1 (\rho(\tau)-\bar{\rho})(\rho(\tau)-\bar{\rho})^\top \,\mathrm{d}\tau\\
        &= \int_0^1 \rho(\tau)\rho(\tau)^\top \,\mathrm{d}\tau - \bar{\rho}\bar{\rho}^\top = Q,
    \end{align*}
    where the last equality follows from the definition of $Q$. From the assumptions of the Lemma, we have that $Q \succeq \alpha I$ for some $\alpha>0$, which implies that $Q$ is invertible and has a unique and invertible square root. Now observe that the function defined by
    $v(\tau)=Q^{-\frac{1}{2}}\tilde{v}(\tau)$, 
    is, by construction, a UPE function. \hfill $\blacksquare$
\end{pf}

For any UPE function $v(\cdot):[0,1]\rightarrow\mathbb{R}^n$, we define
\begin{align}
    M_v:=\sup_{\tau\in[0,1]}\|v(\tau)\|,
\end{align}
which, due to the continuity of $v(\cdot)$ and the compactness of $[0,1]$, is well-defined and finite. }

\section{Algorithm Description: Static Maps}\label{sec:PRLIs}
{ In this section, we introduce our framework in the context of a static cost function. This allows us to explain the core ideas of our framework without the added complexity of the interconnection with the plant. We return to the case of a dynamic cost in Section \ref{sec:dynamic-plants}}.

Consider the following optimization problem:
\begin{align}\label{eq:minimization-problem}
    \min_{u\in \mathcal{U}} \phi(u),
\end{align}
where $\phi:\mathbb{R}^n\to\mathbb{R}$ is a smooth cost function, and  $\mathcal{U}\subseteq\mathbb{R}^n$ is the \emph{feasible} set. If $\mathcal{U}=\mathbb{R}^n$, the simplest iterative algorithm to solve \eqref{eq:minimization-problem} is the discrete-time gradient descent algorithm given by
\begin{align}\label{eq:discrete-time-gd}
    u^+&= u - \gamma \nabla \phi(u),
\end{align}
which, under appropriate assumptions on the function $\phi$ and the \emph{learning rate} $\gamma$, is guaranteed to converge to the set of minimizes of $\phi$ \cite{polyak1987introduction}. In general, depending on the specifications on $\phi$ and $\mathcal{U}$, one may consider a more sophisticated iterative scheme with an augmented state $\eta=(u,w)\in\mathbb{R}^n\times\mathbb{R}^m$, and an update rule of the form:
\begin{align}\label{eq:generic-form-iterative-scheme}
    \eta^+&= g(\eta,\nabla\phi\circ h(\eta)),
\end{align}
for some continuous maps $g$ and $h$. { Examples of such optimization schemes will be studied in Section \ref{sec:stability_properties}.}

\subsection{PRLIs as Gradient Oracles: The Main Idea}
{ As is clear from \eqref{eq:generic-form-iterative-scheme}, all schemes of such form  require the availability of the gradient $\nabla \phi$ at each update time in order to generate the next value of the input $u$. To dispense with this requirement, we now introduce a novel oracle, which is a key component of our framework, for gradient estimation using online measurements of the cost $\phi(u)$.} Namely, we define the map
\begin{align}\label{eq:gradient-estimator}
    \bar{\Phi}_{a}(u)&:=a^{-1}\int_0^1 \phi(u + a v(\tau)) \,v(\tau) \,\mathrm{d}\tau,
\end{align}
where $a$ is a parameter, and $v(\cdot):[0,1]\rightarrow\mathbb{R}^n$ is any UPE function (recall Definition \ref{defn:UPE}). We then have the following proposition.
\begin{prop}[The ``Oracle'' Property]\label{prop:gradient-estimator-properties}
    Let $\phi:\mathbb{R}^n\rightarrow\mathbb{R}$ be continuously differentiable. Then, for any UPE function $v(\cdot)$, the function $(u,a)\mapsto \bar{\Phi}_{a}(u)$ is jointly continuous in $(u,a)\in\mathbb{R}^n\times\mathbb{R}$, and satisfies
    \begin{align}\label{eq:gradient-estimator-simplified}
        \bar{\Phi}_{a}(u) = \nabla \phi(u) + R(u,a),
    \end{align}
    where the error term $R$ satisfies
    \begin{align}\label{eq:unbiased-gradient-estimator-property}
        R(u,0) = 0,~~~~~\forall~u\in\mathbb{R}^n.
    \end{align}
    If, in addition, there exists a constant $L_\phi>0$ such that
    \begin{align}\label{eq:lipschitz-condition}
        \|\nabla \phi(u_1) -\nabla \phi(u_2)\|\leq L_\phi\|u_1-u_2\|,
    \end{align}
    for all $u_1,u_2\in\mathbb{R}^n$, then
    \begin{align}\label{eq:uniform-upper-bound-remainder}
            \lVert R(u,a)\lVert \leq M_v^3 L_\phi |a|,
    \end{align}
    for all $(u,a)\in\mathbb{R}^n\times\mathbb{R}$.
\end{prop}
\begin{pf}
    From Hadamard's Lemma \cite[Lemma 2.8]{nestruev2003smooth},
    \begin{align*}
        \phi(u +a &v(\tau) )=\phi(u)\\
        &+ a  \int_0^1v(\tau)^\top\nabla\phi(u+a\lambda v(\tau))\mathrm{d}\lambda.
    \end{align*}
    Substituting into $\bar{\Phi}_{a}$ and invoking \eqref{eq:unbiased-exploration}, we obtain that
    \begin{align*}
        \bar{\Phi}_{a}(u)&= \nabla\phi(u) + R(u,a),
    \end{align*}
    where $R$ is the vector-valued map given by
    \begin{align*}
        R(u,a)=\!\!\int_0^1\!\!\!\int_0^1 \!\!\! Q(\tau)(\nabla \phi(u+ a\lambda v(\tau))-\nabla\phi(u))\,\mathrm{d}\lambda\,\mathrm{d}\tau,
    \end{align*}
    and $Q(\tau) := v(\tau)v(\tau)^\top$, which proves the first claim of the proposition. To prove the second claim, we observe that, utilizing the sub-multiplicative property of the norm $\|\cdot\|$, we have that
    \begin{align*}
        \lVert Q(\tau)(\nabla &\phi(\bar{u}+a\lambda v(\tau))-\nabla\phi(\bar{u}))\lVert\\
        &\leq \lVert v(\tau)\lVert^2 \lVert \nabla \phi(\bar{u}+2 a\lambda v(\tau))-\nabla\phi(\bar{u})\lVert.
    \end{align*}
    On the other hand, from \eqref{eq:lipschitz-condition}, we obtain that
    \begin{align*}
        \lVert\nabla \phi(\bar{u}+a\lambda v(\tau))-\nabla\phi(\bar{u})\lVert\leq L_\phi |a \lambda|\lVert v(\tau)\lVert,
    \end{align*}
    which implies that
    \begin{align*}
        \|R(u,a)\|\leq \int_0^1\!\! L_\phi |a|\lambda \,\mathrm{d}\lambda\int_0^1\!\!\lVert v(\tau)\lVert^3\mathrm{d}\tau.
    \end{align*}
    By definition of the constant $M_v$, we have that
    \begin{align*}
        \sup_{\tau\in[0,1]}\|v(\tau)\|\leq M_v,
    \end{align*}
    from which the second claim is evident. \hfill $\blacksquare$
\end{pf}
\begin{rem}\normalfont
    Proposition \ref{prop:gradient-estimator-properties} establishes that, for any UPE function $v(\cdot)$, the vector-valued function \eqref{eq:gradient-estimator} is an \emph{unbiased gradient estimator} for $a=0$. Moreover, under a reasonable assumption on the rate of growth of the function $\phi$, i.e. condition \eqref{eq:lipschitz-condition}, the proposition also establishes an upper bound on the remainder term $R$ in terms of $a$ that \emph{is uniform with respect to $u\in\mathbb{R}^n$}. This is a key distinction compared to approximation-based methods such as those based on averaging, which usually only lead to ``semi-global'' gradient approximations that are valid only on compact sets and for ``sufficiently small'' tunable parameters. Instead, property \eqref{eq:gradient-estimator-simplified} is exact (does not rely on asymptotic approximations), { which will allow us to establish global stability under natural assumptions on the cost $\phi$}. \hfill $\square$
\end{rem}

Examination of \eqref{eq:gradient-estimator} shows that $\bar{\Phi}_{a}(u)$ can be written as
\begin{align}
    \bar{\Phi}_{a}(u) = a^{-1} p(1),
\end{align}
for each $(u,a)\in\mathbb{R}^n\times\mathbb{R}_{+}$, where $p:[0,1]\rightarrow\mathbb{R}^n$ is the unique solution of the initial value problem
\begin{align}\label{eq:oracle-definition}
    \frac{dp}{d\tau}&= {\Phi}_{a}(u,\tau) := \phi(u +a v(\tau)) \,v(\tau), & p(0)&=0.
\end{align}
That is, a scaled version of $\bar{\Phi}_{a}$ can be obtained dynamically and in real-time \emph{using an integrator that is reset to zero whenever $\tau=1$}. 
\begin{figure*}[t!]
    \centering
    \begin{subfigure}[b]{0.5\textwidth}
        \centering
        \includegraphics[width=\linewidth]{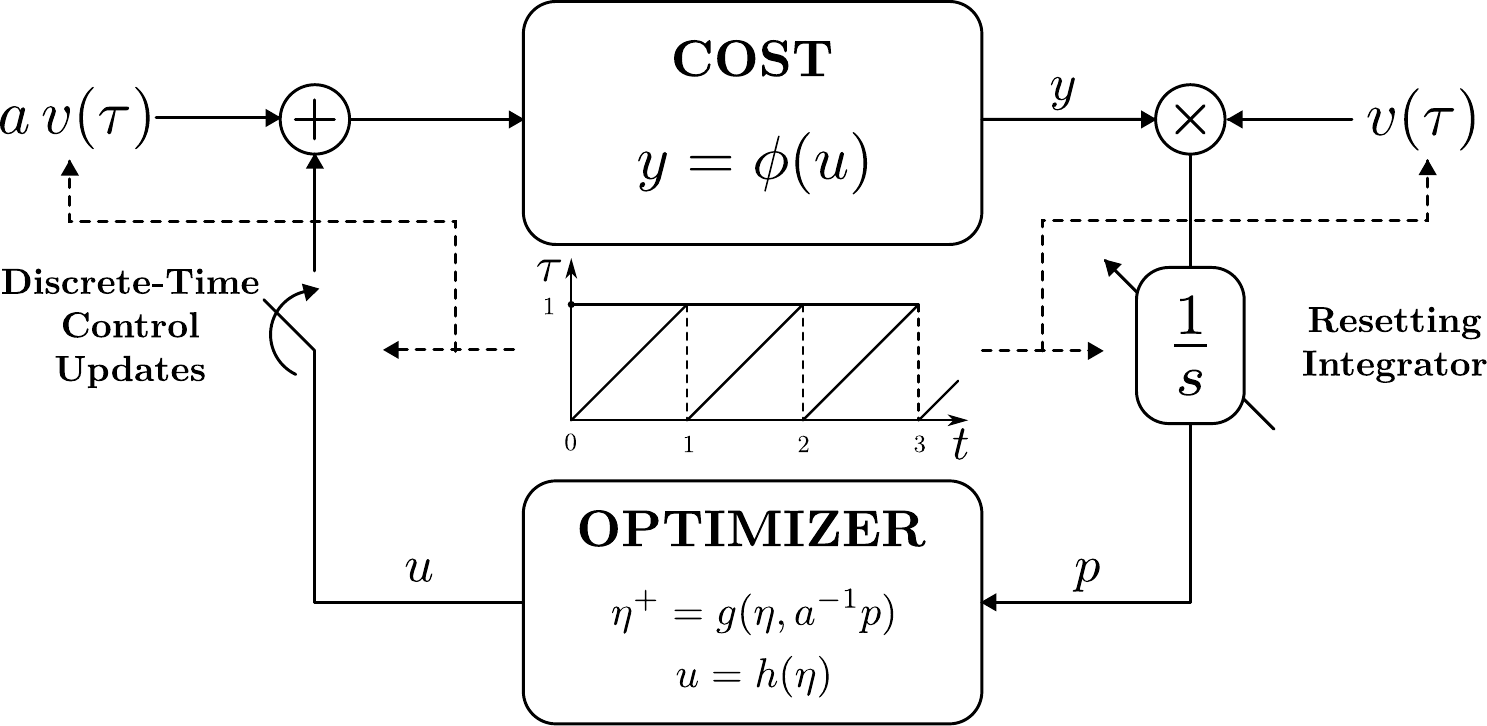}
        \caption{}
        \label{fig:block-diagram}
    \end{subfigure}
    \hfill
    \begin{subfigure}[b]{0.24\textwidth}
        \centering
    \includegraphics[width=\linewidth]{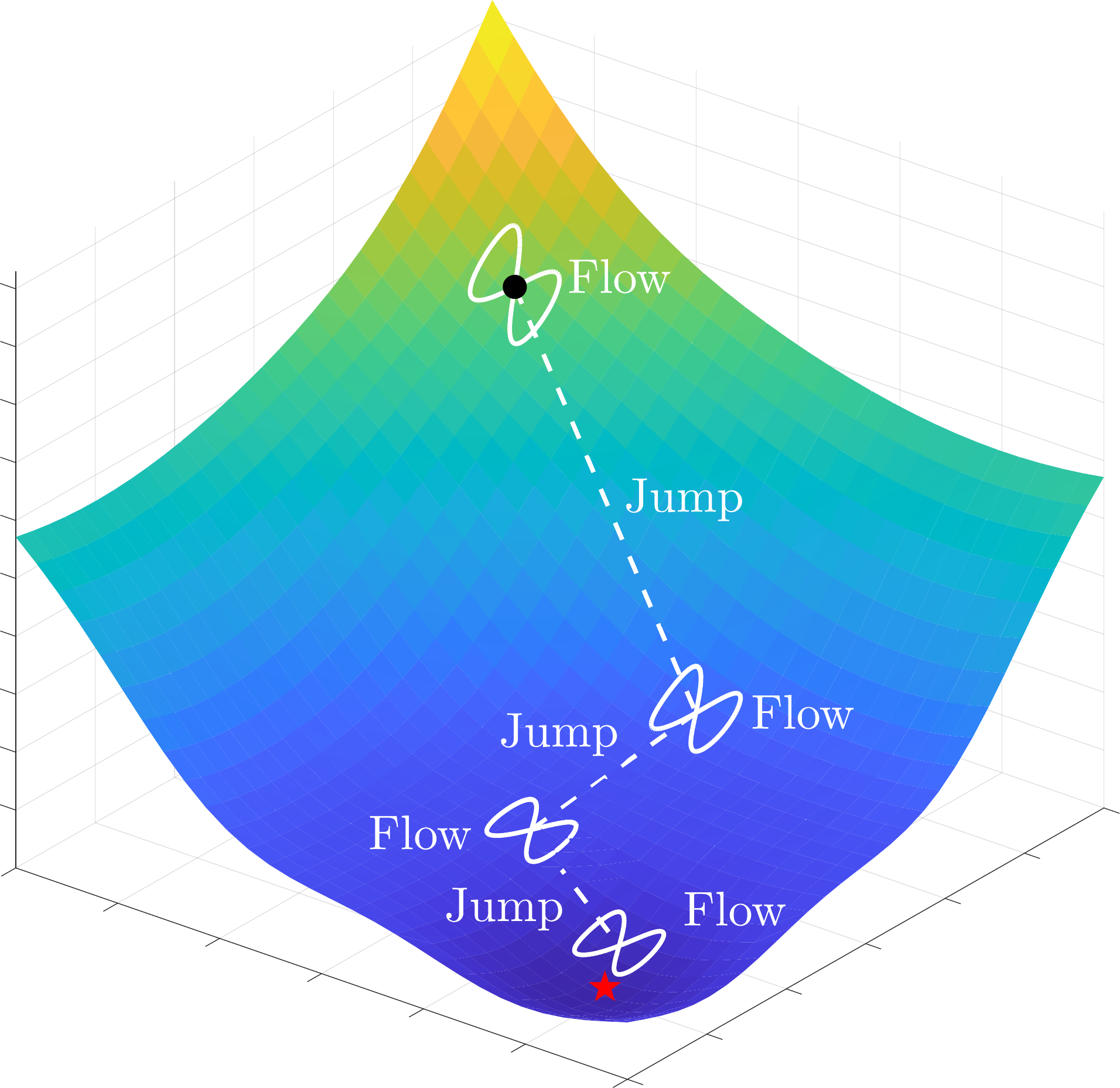}
    \caption{}
    \label{fig:function-figure-3d}
    \end{subfigure}
    \hfill
    \begin{subfigure}[b]{0.24\textwidth}
        \centering
    \includegraphics[width=\linewidth]{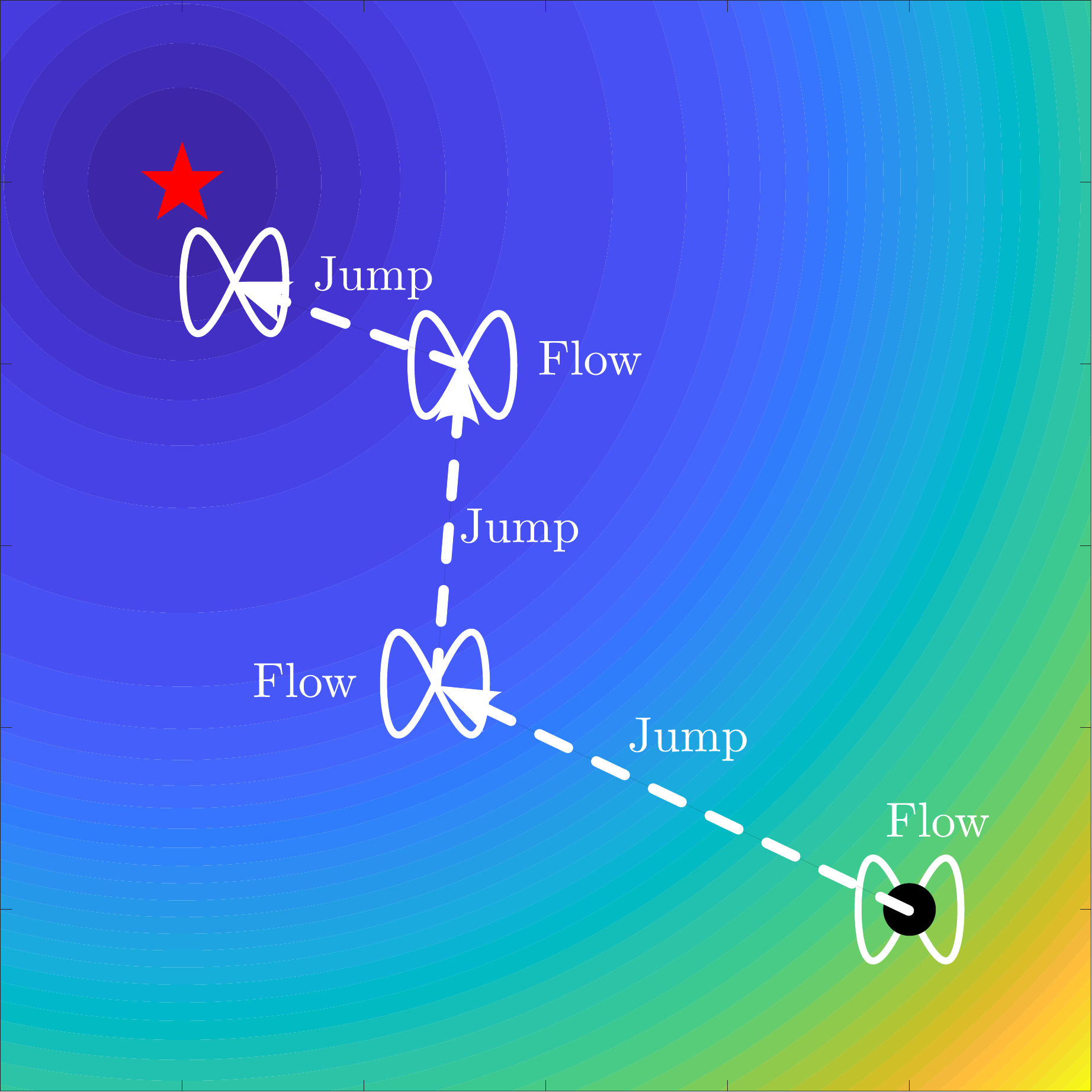}
    \caption{}
    \label{fig:function-figure}
    \end{subfigure}
    \caption{(a): A block diagram description of the interconnection between the proposed model-free optimization algorithm and a static cost function. (b): a 3D sketch of the behavior of the optimization parameter $u$ and the cost $\phi(u)$ under the HDS $\mathcal{H}$. (c): a 2D sketch of the typical behavior of the optimization parameter $u$ under the HDS $\mathcal{H}$.}
\end{figure*}
\subsection{Formulation as Hybrid Systems}
The previous observations naturally suggest a family of model-free optimization algorithms that can be formulated as the HDS $\mathcal{H}=(C,F,D,G)$ with state $\xi=(\eta,p,\tau)$, where $\eta$ is the main state of the optimization algorithm, $p$ is the state of the integrator, and $\tau$ is a resetting timer. The dynamics of the system are then characterized by the following flow and jump sets:
\begin{subequations}\label{eq_SCRLI}
\begin{align}
    C&= \mathbb{R}^{n+m}\times\mathbb{R}^n\times[0,1], \\
    D&= \mathbb{R}^{n+m}\times\mathbb{R}^n\times\{1\},
\end{align}
and the following flow and jump maps:
\begin{align}
    F(\xi)&=(0,{\Phi}_{a}(h(\eta),\tau),1), \\
    G(\xi)&=(g(\eta,a^{-1}p), 0, 0).
\end{align}
\end{subequations}
We note that the HDS $\mathcal{H}$ satisfies Assumption \ref{asmp:regularity}, i.e. the HDS $\mathcal{H}$ is a well-posed HDS. A high-level block diagram description of the algorithms is shown in Figure \ref{fig:block-diagram}, which emphasizes the role of the integrator with state $p$ and the resetting timer $\tau$, which coordinates the resets of the integrator, the function $v$, and the updates of the optimization algorithm. { We also provide a sketch in figures \ref{fig:function-figure} and \ref{fig:function-figure-3d} depicting a typical trajectory of the optimization parameter $u$ under the evolution of the HDS $\mathcal{H}$. In figures \ref{fig:function-figure} and \ref{fig:function-figure-3d}, the red star represents the minimizer, the black circle represents the initial condition, the continuous lines represent the portions of the trajectory corresponding to flows (performing exploration), and the dashed lines represent the jumps (performing exploitation).}

To study the qualitative properties of the hybrid system $\mathcal{H}$, we define a ``target" HDS $\bar{\mathcal{H}}=(\bar{C},\bar{F},\bar{D},\bar{G})$ with state $\bar{\xi}=(\bar{\eta},\bar{\tau})$, flow and jump sets given by
\begin{align*}
    \bar{C}&= \mathbb{R}^{n+m}\times[0,1], &
    \bar{D}&= \mathbb{R}^{n+m}\times\{1\},
\end{align*}
and flow and jump maps given by
\begin{align*}
    \bar{F}(\bar{\xi})&=(0, 1), &
    \bar{G}(\bar{\xi})&=(g(\bar{\eta},\bar{\Phi}_{a}(h(\bar{\eta}))),0),
\end{align*}
where $\bar{\Phi}_{a}$ is the oracle defined by \eqref{eq:gradient-estimator}. Then, we have the following Theorem.
\begin{thm}\label{thm:solution_properties}\normalfont
    Let $\phi:\mathbb{R}^n\rightarrow\mathbb{R}$ be continuously differentiable, and let $g(\cdot)$ and $h(\cdot)$ be continuous. Let $a\in\mathbb{R}_{+}$ and suppose that $\xi=(\eta,p,\tau)$ is a maximal solution of $\mathcal{H}$ with $\xi(0,0)\in C\cup D$. Then, the following holds:
    \begin{enumerate}
        \item The solution $\xi$ is complete, i.e.
        \begin{align}\label{eq:maximal-solution-domain}
            \text{dom}(\xi)=\bigcup_{j=0}^\infty ([t_j,t_{j+1}],j),
        \end{align}
        with $t_1\leq 1$, and $t_{j+1}-t_{j}=1$ for all $j\geq 1$;
        \item Any solution satisfies
        \begin{align}\label{eq:maximal-solution-aux-state-properties}
            \|p(t,j)\|\leq M_v \max_{\tilde{u}\in a M_v\mathbb{B}}|\phi(h(\eta(t_j,j))+ \tilde{u})|,
        \end{align}
        for all $(t,j)\in\text{dom}(\xi)$ with $j\geq 1$;
        \item The hybrid arc $\bar{\xi}=(\bar{\eta},\bar{\tau})$, where
        \begin{subequations}\label{eq:maximal-solution-properties}
            \begin{align}
                \bar{\eta}(t-t_1,j-1)&= \eta(t,j), \\
                \bar{\tau}(t-t_1,j-1)&= \tau(t,j), 
            \end{align}
        \end{subequations}
        for all $(t,j)\in\text{dom}(\xi)$ with $j\geq 1$, is a solution of $\bar{\mathcal{H}}$.
    \end{enumerate}
\end{thm}
\begin{pf}
    We observe that $C\cap D = D$. Let $\xi(0,0)\in C\cup D$. Then, exactly one of the following cases is true: (C1) $\xi(0,0)\in D$; or (C2) $\xi(0,0)\in C \backslash D$. If (C1) is true then, by construction of the set $C$, a maximal solution must jump since in that case $F(\xi(0,0))\not \in T_{\xi(0,0)} C$, where $T_{\xi} C$ is the tangent cone to $C$ at the point $\xi$. Thus, $t_1=0$ and, from the jump map, we obtain that $p(t_1,1) = 0$, $\tau(t_1,1)= 0$, and by continuity of $\phi$, $g$, and $h$,
    \begin{align*}
        \lVert \eta(t_1,1)\lVert &< +\infty.
    \end{align*}
    On the other hand, if (C2) is true, then a maximal solution must flow since $F(\xi(0,0))\in T_{\xi(0,0)} C$. The flow map of $\mathcal{H}$ is explicitly integrable. Through direct integration, we see that $u(t,0)= u(0,0)$, $w(t,0)= w(0,0)$, and 
    \begin{align*}
        \tau(t,0)&= \tau(0,0) + t,\\
        p(t,0)&= p(0,0) +  {\int_{0}^{t}}\dot{p}(s,0)\,\mathrm{d}s, 
    \end{align*}
    for all $(t,j)\in\text{dom}(\xi)$, where $\dot{p}(t,0)$ is
    \begin{align*}
        \dot{p}(t,0)&= {\Phi}_{a}(h(\eta(0,0)),\tau(t,0)).
    \end{align*}
    Explicit computation shows that
    \begin{align*}
        \|\dot{p}(t,0)\|&\leq |\phi(h(\eta(0,0)) +2 a v(\tau(t,0)))| \|v(\tau(t,0))\|\\
                        &\leq M_v \max_{\tau\in[0,1]}|\phi(h(\eta(0,0))+2 a v(\tau))|=: B_0,
    \end{align*}
    Hence, we have that
    $$\lVert p(t,0)\lVert\leq \lVert p(0,0)\lVert + B_0 t,$$ 
    for all $t\geq 0$ such that $(t,0)\in\text{dom}(\xi)$.
    On the other hand, since $\xi(0,0)\in C\backslash D$, it follows that $\tau(0,0)\in[0,1)$, and so, from the definition of the flow map,
    \begin{align*}
        \tau(1-\tau(0,0),0)&= \tau(0,0) + 1-\tau(0,0) = 1.
    \end{align*}
    Therefore, $\xi(1-\tau(0,0),0)\in D$, which implies that $t_1=1-\tau(0,0) \leq 1$. It follows that
    $$\lVert p(t,0)\lVert\leq \lVert p(0,0)\lVert + B_0,$$ 
    and, by substituting into the jump map, we obtain that $p(t_1,1)=0$, $ \tau(t_1,1)=0$, and 
    \begin{align*}
         \lVert \eta(t_1,1)\lVert &= \lVert g(\eta(t_1,0),a^{-1} p(t_1,0))\lVert < +\infty.
    \end{align*}
    To summarize, we showed that in both of the cases (C1) and (C2), we have that $p(t_1,1)=0$, $\tau(t_1,1)=0$, and $\lVert \eta(t_1,1)\lVert < +\infty$, which implies that $\xi(t_1,1)\in C\backslash D$. By repeating the argument we used for the case (C2), we can show that $\xi$ is complete, that its domain is of the form \eqref{eq:maximal-solution-domain}, that 
    $$t_{j+1}-t_j=1,\quad p(t_j,j)=0, \quad\tau(t_j,j)=0,$$
    for all $j\geq 1$, and that \eqref{eq:maximal-solution-aux-state-properties} is satisfied, $\forall (t,j)\in\text{dom}(\xi)$ with $j\geq 1$. This proves the first and second claims of the proposition. To prove the third claim, we define the hybrid arc $\bar{\xi}=(\bar{\eta},\bar{\tau})$, where
    \begin{align*}
        \bar{\eta}(t-t_1,j-1)&=\eta(t,j), \\ \bar{\tau}(t-t_1,j-1)&= \tau(t,j), 
    \end{align*}
    $\forall (t,j)\in\text{dom}(\xi)$ with $j\geq 1$. Explicit computation gives
    \begin{align*}
        \text{dom}(\bar{\xi}\,)&=\bigcup_{j=0}^\infty ([t_j,t_{j+1}],j),\\
        \bar{\eta}(t_{j+1},j+1)&=g(\bar{\eta}(t_{j},j),\bar{\Phi}_{a}(h(\bar{\eta}(t_{j},j)))),
    \end{align*}
    for all $\forall (t,j)\in\text{dom}(\bar{\xi})$, where $ t_j= j$, $\bar{\tau}(t_j,j)=0$, $\bar{\eta}(t,j)=\bar{\eta}(t_j,j)$, and $\bar{\Phi}_{a}$ is the oracle defined in \eqref{eq:gradient-estimator}, see also \eqref{eq:gradient-estimator-simplified}. It follows from the definition of $\bar{\mathcal{H}}$ that $\bar{\xi}$ is a complete solution of $\bar{\mathcal{H}}$ starting from the initial condition $\bar{\xi}(0,0)=(\bar{u}(0,0),0)\in \bar{C}\backslash \bar{D}$. \hfill $\blacksquare$
\end{pf}

\vspace{-0.3cm}
Thanks to Theorem \ref{thm:solution_properties}, the properties of the solutions of the HDS $\bar{\mathcal{H}}$ can be used to inform the analysis of the properties of the HDS $\mathcal{H}$. Indeed, the following Corollary is a direct consequence of the Theorem.
\begin{cor}\normalfont\label{cor:corollary-to-solution-properties}
    Suppose that, for some $a\in\mathbb{R}_{+}$, there exists a class $\mathcal{KL}$ function $\beta$ and a closed set $\bar{\mathcal{A}}\subset \bar{C}\cup \bar{D}$ such that every maximal solution of the HDS $\bar{\mathcal{H}}$ satisfies the uniform $\mathcal{KL}$ bound
    \begin{align}
        \|\bar{\xi}(t,j)\|_{\bar{\mathcal{A}}}< \beta(\|\bar{\xi}(0,0)\|_{\bar{\mathcal{A}}},t+j),
    \end{align}
    for all $(t,j)\in\mathrm{dom}(\bar{\xi})$. Then, every maximal solution $\xi=(\eta,\tau,p)$ of $\mathcal{H}$ is complete and, with $\zeta=(\eta,\tau)$, satisfies
    \begin{align}
        \|\zeta(t,j)\|_{\bar{\mathcal{A}}}< \beta(\|\zeta(t_1,1)\|_{\bar{\mathcal{A}}},t+j-t_1-1),
    \end{align}
    for all $(t,j)\in\mathrm{dom}(\xi)$ with $j\geq 1$.
\end{cor}

The following result establishes the existence of a compact UGAS stable set for the HDS $\mathcal{H}$ whenever the HDS $\bar{\mathcal{H}}$ has a UGAS compact subset for some $a\in\mathbb{R}_{+}$. 

\begin{thm}\normalfont\label{thm:orig-sys-is-ugas}
    Let $\phi(\cdot)$ be continuously differentiable, and let $g(\cdot)$ and $h(\cdot)$ be continuous. Suppose that, for some $a\in\mathbb{R}_+$, there exists a compact subset $\bar{\mathcal{A}}\subset\mathbb{R}^{n+m}$ such that $\bar{\mathcal{A}}\times[0,1]$ is UGAS for the HDS $\bar{\mathcal{H}}$. Then there exists a compact subset $\mathcal{A}\subset \mathbb{R}^{n+m}\times\mathbb{R}^{n}$ such that 
    \begin{align*}
    \mathcal{A}\subseteq \bar{\mathcal{A}}\times\{\|p\|\leq \max_{\eta\in\bar{\mathcal{A}}}\max_{\tilde{u}\in M_v \mathbb{B}} M_v|\phi(h(\eta)+a\tilde{u})|\},
    \end{align*}
    and $\mathcal{A}\times[0,1]$ is UGAS for the HDS $\mathcal{H}$.
\end{thm}
{ Before we proceed with the proof of Theorem \ref{thm:orig-sys-is-ugas}, we highlight some of its important implications. Specifically, Theorem \ref{thm:orig-sys-is-ugas} provides a \emph{uniform global} stability result for a compact set $\mathcal{A}$. Moreover, it provides an explicit characterization of a bounding box that contains this compact set. We emphasize that this is in contrast to existing literature where, at best, semi-global practical stability can be established.}

\begin{figure}[t!]
    \centering
    \includegraphics[width=0.75\linewidth]{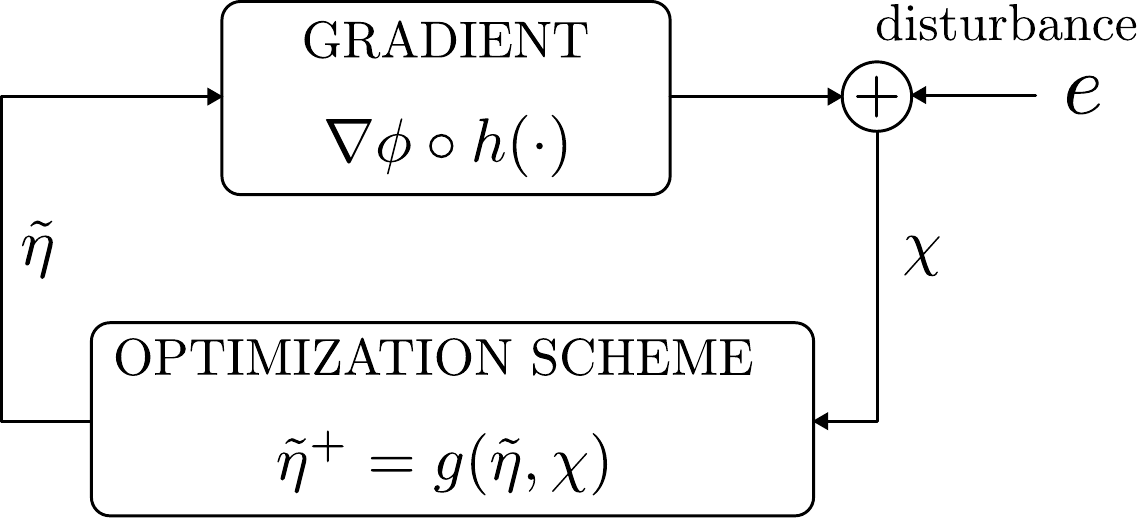}
    \caption{A block diagram of the discrete-time system \eqref{eq:discrete-time-system-with-input}.}
    \label{fig:discrete-time-block-diagram}
\end{figure}

\vspace{-0.3cm}
\begin{pf}
    Since $\phi(\cdot)$ is continuously differentiable and $g(\cdot)$ and $h(\cdot)$ are continuous, it follows that the flow and jump maps of the HDS $\mathcal{H}$ and the HDS $\bar{\mathcal{H}}$ are continuous and, therefore, both systems are well-posed. From Theorem \ref{thm:solution_properties}, we know that every maximal solution $\xi = (\eta,p,\tau)$ of the HDS $\mathcal{H}$ is complete, has a domain of the form \eqref{eq:maximal-solution-domain}, and is such that the hybrid arc defined by \eqref{eq:maximal-solution-properties} is a solution of the HDS $\bar{\mathcal{H}}$. Since, by assumption, there exists a compact set $\bar{\mathcal{A}}\subset\mathbb{R}^{n+m}$ such that the set $\bar{\mathcal{A}}\times[0,1]$ is UGAS for the HDS $\bar{\mathcal{H}}$, we know from Corollary \ref{cor:corollary-to-solution-properties} that there exists a class $\mathcal{KL}$ function $\beta$ such that every solution $\xi$ of the HDS $\mathcal{H}$ satisfies
\begin{align}\label{eq:eta_ultimate_bound}
    \|\eta(t,j)\|_{\bar{\mathcal{A}}}< \beta(\|\eta(t_1,1)\|_{\bar{\mathcal{A}}},t+j-t_1-1),
\end{align}
for all $(t,j)\in\mathrm{dom}(\xi)$ with $j\geq 1$. On the other hand, from Theorem \ref{thm:solution_properties}, we know that every solution $\xi = (\eta,p,\tau)$ of the HDS $\mathcal{H}$ also satisfies \eqref{eq:maximal-solution-aux-state-properties}, i.e. that
\begin{align}\label{eq:p_ultimate_bound}
    \|p(t,j)\|\leq M_v \max_{\tilde{u}\in a M_v\mathbb{B}}|\phi(h(\eta(t_j,j))+ \tilde{u})|,
\end{align}
for all $(t,j)\in\mathrm{dom}(\xi)$ with $j\geq 1$. Combining \eqref{eq:eta_ultimate_bound} and \eqref{eq:p_ultimate_bound}, we see that the $\omega$-limit set $\Omega(C\cup D)$ of the HDS $\mathcal{H}$ is compact and is contained within the compact set $\mathcal{K}\times[0,1]$, where $\mathcal{K}$ is defined by
\begin{align*}
    \mathcal{K}:= \bar{\mathcal{A}}\times\{\|p\|\leq M_v \max_{\eta\in\bar{\mathcal{A}}}\max_{\tilde{u}\in \mathbb{B}} M_v|\phi(h(\eta)+a\tilde{u})|\}.
\end{align*}
That is, since every maximal solution of $\mathcal{H}$ is complete, we have that
\begin{align*}
    \lim_{t+j\rightarrow +\infty} |(\eta(t,j),p(t,j))|_{\mathcal{K}} = 0.
\end{align*}
Therefore, from the definition of uniform global recurrence in Definition \ref{defn:recurrence}, we see that if $\mathcal{O}\subset\mathbb{R}^{n+m}\times\mathbb{R}^n$ is any bounded and open set containing the compact set $\mathcal{K}$, then the subset $\mathcal{O}\times[0,1]$ is uniformly globally recurrent for the HDS $\mathcal{H}$. However, since $\mathcal{H}$ is well-posed, it follows that \cite[Proposition 4]{subbaraman2016equivalence} there exists a compact set $\mathcal{A}$ such that $\mathcal{A}\times[0,1]$ is UGAS for the HDS $\mathcal{H}$. Since every maximal solution of the HDS $\mathcal{H}$ is complete, it follows that $\mathcal{A}\times[0,1]$ is necessarily a subset of the $\omega$-limit set $\Omega(C\cup D)$, i.e. $\mathcal{A}\subset \mathcal{K}$, which concludes the proof. \hfill $\blacksquare$
\end{pf}
\begin{rem}\normalfont
    Although Theorem \ref{thm:orig-sys-is-ugas} establishes stability of the HDS $\mathcal{H}$ based on the stability of the HDS $\bar{\mathcal{H}}$, no assertions are made concerning the uniformity of the stability properties of the HDS $\mathcal{H}$ with respect to variations in the parameter $a$. The obstruction to such uniformity is the fact that, unlike the jump map $\bar{G}$ of the HDS $\bar{\mathcal{H}}$, which is continuous with respect to $a$, for all $a\in\mathbb{R}$, the jump map $G$ of the HDS $\mathcal{H}$ is not continuous with respect to $a$ at $a=0$. Indeed, from the definition of the flow map $F$,
    \begin{align*}
        p(t_1,0)&= p(0,0) + \int_0^{t_1} \Phi_a(h(\eta(0,0)),\tau(t,0))\,\mathrm{d}t
    \end{align*}
    and so, if $t_1\neq 1$ or if $p(0,0)\neq 0$, then $a^{-1}p(t_1,0)=O(a^{-1})$. Thus, in the limit $a\rightarrow 0^+$, the first jump may push the system further and further away from the set $\bar{\mathcal{A}}$ unless the system is properly initialized. { We note, however, that this issue is unavoidable and is commonly encountered in sampled data and event-triggered seeking systems, see \cite{poveda2017robust}.} On the other hand, because the states $p$ and $\tau$ are internal states of the oracle, it is always possible to initialize the system in a way that guarantees convergence without any prior knowledge of $\phi$. \hfill $\square$
\end{rem}
Thanks to Theorem \ref{thm:orig-sys-is-ugas}, the task of establishing the stability of the HDS $\mathcal{H}$ reduces to proving the existence, for some $a\in\mathbb{R}_+$, of a UGAS compact set $\mathcal{A}$ for the HDS $\bar{\mathcal{H}}$. The advantage of this approach is that the HDS $\bar{\mathcal{H}}$ is essentially a \emph{perturbed} discrete-time iterative scheme. Indeed, if we define the discrete-time system
\begin{align}\label{eq:discrete-time-system}
    \tilde{\eta}^+&= g(\tilde{\eta},\Phi_a\circ h(\tilde{\eta})),
\end{align}
then we have the following Lemma\footnote{The proof of Lemma \ref{lem:solution-relation-discrete-hybrid} is almost identical to part of the proof of Theorem \ref{thm:solution_properties}, and so is omitted. }.
\begin{lem}\normalfont\label{lem:solution-relation-discrete-hybrid}
    Let $\bar{\xi}=(\bar{\eta},\tau)$ be any maximal solution of $\bar{\mathcal{H}}$. Then, $\bar{\xi}$ is complete, its domain is of the form \eqref{eq:maximal-solution-domain}, 
    \begin{align*}
        \bar{\eta}(t,j)&= \bar{\eta}(t_{j},j),
    \end{align*}
    for all $(t,j)\in\mathrm{dom}(\bar{\xi})$, and the function $\tilde{\eta}:\mathbb{N}\rightarrow\mathbb{R}^n\times\mathbb{R}^m$ defined by $\tilde{\eta}(j)= \bar{\eta}(t_{j},j)$
    is a solution of the discrete-time system \eqref{eq:discrete-time-system}. \hfill $\square$
    
\end{lem}
On the other hand, the discrete-time system \eqref{eq:discrete-time-system} can be written in the form
\begin{align*}
    \tilde{\eta}^+&= g(\tilde{\eta},\nabla\phi\circ h(\tilde{\eta})+e(\tilde{\eta})), 
\end{align*}
where $e$ is given by
\begin{align*}
    e(\tilde{\eta})&=\Phi_a\circ h(\tilde{\eta})-\nabla\phi\circ h(\tilde{\eta})
\end{align*}
for all $(\tilde{\eta},a)\in\mathbb{R}^n\times\mathbb{R}^m\times\mathbb{R}$. If $\phi$ satisfies the Lipschitz condition \eqref{eq:lipschitz-condition}, then Proposition \ref{prop:gradient-estimator-properties} implies that the error term $e(\tilde{\eta})$ satisfies the uniform upper bound
\begin{align*}
    \|e(\tilde{\eta})\|\leq M_v^3 L_\phi |a|,
\end{align*}
for all $(\tilde{\eta},a)\in\mathbb{R}^{n+m}\times\mathbb{R}$. The preceding discussion suggests a natural route to establish the existence of a compact set $\bar{\mathcal{A}}$ such that $\bar{\mathcal{A}}\times[0,1]$ is UGAS for the HDS $\bar{\mathcal{H}}$. The first step along this route is to establish \emph{Input-to-State Stability} (ISS) for the system
\begin{align}\label{eq:discrete-time-system-with-input}
    \tilde{\eta}^+&= g(\tilde{\eta},\nabla\phi\circ h(\tilde{\eta})+e),
\end{align}
wherein $e \in\mathbb{R}^n$ is treated as an input (see Figure \ref{fig:discrete-time-block-diagram} for a block diagram description of \eqref{eq:discrete-time-system-with-input}). Recall that the system \eqref{eq:discrete-time-system-with-input} is said to be ISS if there exists a compact set $\mathcal{A}_d\subset\mathbb{R}^{n+m}$, a class $\mathcal{KL}$ function $\beta$, and a class $\mathcal{K}_\infty$ function $\alpha$ such that every solution of \eqref{eq:discrete-time-system-with-input} satisfies
\begin{align*}
    |\tilde{\eta}(j)|_{\mathcal{A}_d}\leq \beta(|\tilde{\eta}(0)|_{\mathcal{A}_d},j) + \alpha(\|e\|_\infty).
\end{align*}
If such a property holds and the gradient $\nabla\phi$ satisfies a Lipschitz condition of the form \eqref{eq:lipschitz-condition}, it follows that every solution of the discrete-time system \eqref{eq:discrete-time-system} satisfies
\begin{align*}
    |\tilde{\eta}(j)|_{\mathcal{A}_d}\leq \beta(|\tilde{\eta}(0)|_{\mathcal{A}_d},j) + \alpha(M_v^3 L_\phi|a|),
\end{align*}
for all $a\in\mathbb{R}$. By invoking Lemma \ref{lem:solution-relation-discrete-hybrid}, we obtain that the HDS $\bar{\mathcal{H}}$ is $\kappa$-UGUB with respect to the compact set $\mathcal{A}_d\times[0,1]$, where $\kappa = \alpha(M_v^3 L_\phi|a|)$,
is the ultimate bound. The preceding argument forms the main part of the proof for the following theorem.
\begin{thm}\label{thm:avg-sys-is-ugas}\normalfont
    Let $\phi$ be a continuously differentiable function that satisfies \eqref{eq:lipschitz-condition}, and let $g$ and $h$ be continuous. Suppose that the discrete-time system \eqref{eq:discrete-time-system-with-input} is ISS with respect to a compact subset $\mathcal{A}_d$. Then, for every $a\in\mathbb{R}^+$, there exists a compact subset $\bar{\mathcal{A}}$ such that
    \begin{align}
        \bar{\mathcal{A}}\subset \mathcal{A}_d+\alpha(M_v^3 L_\phi|a|)\mathbb{B},
    \end{align}
    and $\bar{\mathcal{A}}\times[0,1]$ is UGAS for the HDS $\bar{\mathcal{H}}$. 
\end{thm}
\begin{pf}
    From the discussion that preceded the statement of the Theorem, we know that the HDS $\bar{\mathcal{H}}$ is $\alpha(M_v^3 L_\phi|a|)$-UGUB with respect to the compact set $\mathcal{A}_d\times [0,1]$, and so every solution $\bar{\xi}=(\bar{\eta},\bar{\tau})$ of the HDS $\bar{\mathcal{H}}$ satisfies 
\begin{align*}
    |\bar{\eta}(t,j)|_{\mathcal{A}_d}\leq \beta(|\bar{\eta}(0,0)|_{\mathcal{A}_d},t+j) + \alpha(M_v^3 L_\phi|a|).
\end{align*}
Since the HDS $\bar{\mathcal{H}}$ is a well-posed HDS, it is UGUB if and only if there exists a compact set that is UGAS \cite[Propositions 2 \& 4]{subbaraman2016equivalence}. Due to the structure of the flow and jump maps, namely that
\begin{align*}
    \begin{cases}
        \bar{\tau}\in [0,1] & \dot{\bar{\tau}}=1,\\
        \bar{\tau}\in\{1\} & \bar{\tau}^+=0,
    \end{cases}
\end{align*}
it is clear that such a compact set is of the form $\bar{\mathcal{A}}\times[0,1]$ for some compact set $\bar{\mathcal{A}}\subset\mathbb{R}^{n\times m}$. On the other hand, we see that the $\omega$-limit set $\Omega(\bar{C}\cup\bar{D})$ for the HDS $\bar{\mathcal{H}}$ is necessarily contained within the set
\begin{align*}
    \left(\mathcal{A}_d+\alpha(M_v^3 L_\phi|a|) \mathbb{B} \right)\times[0,1].
\end{align*}
We conclude the proof by noting that if the set $\bar{\mathcal{A}}\times[0,1]$ is UGAS for the HDS $\bar{\mathcal{H}}$, and since every solution of $\bar{\mathcal{H}}$ is complete, then $\bar{\mathcal{A}}\times[0,1]$ is necessarily a subset of the $\omega$-limit set $\Omega(\bar{C}\cup\bar{D})$, i.e. $\bar{\mathcal{A}}\subset\mathcal{A}_d + \alpha(M_v^3 L_\phi|a|) \mathbb{B}$. \hfill $\blacksquare$
\end{pf}

\vspace{-0.2cm}
The compendium of results established so far hinges on the ISS properties of the discrete-time system \eqref{eq:discrete-time-system-with-input}. Nevertheless, recent advances in the literature on control-theoretic stability analysis of optimization algorithms provide a rich toolkit to establish such conditions \cite{lessard2022analysis,hauswirth2024optimization}. In the next section, we show that the ISS property indeed holds for the discrete-time system \eqref{eq:discrete-time-system-with-input} associated with different well-known optimization schemes under natural conditions on the function $\phi$.

\section{Input-to-State Stability of Discrete-Time Optimization Schemes}\label{sec:stability_properties}
In this section, we establish the ISS property of the discrete-time system \eqref{eq:discrete-time-system-with-input} for three common applications of systems of the form \eqref{eq:generic-form-iterative-scheme}: a) gradient descent for \emph{unconstrained} optimization; b) Polyak's Heavy Ball method for \emph{accelerated optimization}; and c) projected gradient descent for \emph{safe optimization}. In this way, it is shown that all these methods can be integrated into the model-free framework considered in Section 3 and the stability results of Theorems \ref{thm:orig-sys-is-ugas} and \ref{thm:avg-sys-is-ugas}.

\subsection{Input-to-State Stability of Gradient Descent}
We begin our investigation with the simplest optimization scheme: gradient descent with a constant step-size, which corresponds to the discrete-time system
\begin{align}\label{eq:discrete-time-gd-with-input}
    \tilde{u}^+&= \tilde{u}-\gamma (\nabla\phi(\tilde{u}) + e).
\end{align}
Throughout the current subsection, we impose the following standing assumption.

\begin{asmp}\label{asmp:cost_regularity}\normalfont The function
    $\phi$ is continuously differentiable, and there exists a constant $L_\phi>0$, a point $u^\star\in\mathbb{R}^n$, and class $\mathcal{K}_\infty$ functions $\mu_1$ and $\mu_2$, such that, for all $u,\tilde{u}\in\mathbb{R}^n$, $\nabla \phi (u)=0$ if and only $u=u^\star$, and
    \begin{align}
    \mu_1 (\lVert u - u^\star\lVert)&\leq \phi(u)-\phi(u^\star),\label{eq:cost_is_pd}\\
    \mu_2(\phi(u)-\phi(u^\star))&\leq \lVert \nabla \phi(u)\lVert, \label{eq:polyak_lojasiewicz}\\
    \lVert\nabla \phi(u)-\nabla\phi(\tilde{u})\lVert&\leq L_\phi\lVert u-\tilde{u}\lVert. \label{eq:lipschitz_condition}
    \end{align}
\end{asmp}
\begin{rem}\normalfont\normalfont
    The inequality \eqref{eq:cost_is_pd} implies that the function $\phi$ is radially unbounded and, up to a constant, is positive definite with respect to $u^\star$. Conversely, any radially unbounded and positive definite function satisfies the inequality \eqref{eq:cost_is_pd} for some class $\mathcal{K}_\infty$ function \cite[Lemma 4.3]{khalil2002nonlinear}. 
    The inequality \eqref{eq:polyak_lojasiewicz} is a generalization of the Polyak-Lojasiewicz (PL) inequality and is a standard assumption in gradient-based optimization \cite{sontag2022remarks}. Finally, the Lipschitz condition \eqref{eq:lipschitz_condition} is indispensable to guarantee global convergence of gradient-based optimization algorithms. Indeed, without the Lipschitz condition \eqref{eq:lipschitz_condition}, there are no guarantees on the convergence of discrete-time gradient descent \cite[p.22]{polyak1987introduction}.
\end{rem}
Under Assumption \ref{asmp:cost_regularity}, we have the following Proposition.
\begin{prop}\label{prop:gd-is-iss}
    Let Assumption \ref{asmp:cost_regularity} be satisfied and let $0<\gamma L_\phi<2$. Then, the discrete time system \eqref{eq:discrete-time-gd-with-input} is ISS with respect to $\mathcal{A}_d=\{u^\star\}$. 
\end{prop}
\begin{pf}
    We define the Lyapunov function candidate
    \begin{align}\label{eq:lyapunov-function}
        V(\tilde{u}):=\phi(\tilde{u})-\phi(u^\star).
    \end{align}
    Using the notation $\Delta V:=V(\tilde{u}^+)- V(\tilde{u})$, we obtain from the Descent Lemma \cite[Proposition 6.1.2]{bertsekas2015convex} that
    \begin{align}\label{eq:lyap-increment-1}
        \Delta V \leq \nabla V(\tilde{u})^\top(\tilde{u}^+-\tilde{u})+\frac{L_\phi}{2}\|\tilde{u}^+-\tilde{u}\|^2,
    \end{align}
    for all $\tilde{u}\in\mathbb{R}^n$. Substituting from \eqref{eq:discrete-time-gd-with-input} into \eqref{eq:lyap-increment-1}, and rearranging terms, we obtain the inequality
    \begin{align*}
        \Delta V &\leq -2 c_2\|\nabla V(\tilde{u})\|^2 + c_1 \|e\| \|\nabla V(\tilde{u})\|+c_0 \|e\|^2,
    \end{align*}
    where the constants $c_0 ,c_1, c_2 $ are given by
    \begin{align*}
        c_2&= \frac{\gamma}{2} (1-\gamma L_\phi/2), & c_1&= \gamma (1+\gamma L_\phi), & c_0&= \frac{\gamma^2 L_\phi}{2}.
    \end{align*}
    It follows that
    \begin{align}\label{eq:lyap-increment-2}
        \Delta V(\tilde{u}) \leq &-c_2\lVert\nabla V(\tilde{u})\lVert^2< 0,
    \end{align}
    for all $\tilde{u}\in\mathbb{R}^n$ that satisfy the condition
    \begin{align}
        \|\nabla V(\tilde{u})\|\geq c_3\|e\|,
    \end{align}
    where $c_3$ is the positive constant given by $c_3=(c_1^2+4 c_0 c_2)^{\frac{1}{2}}+c_1)/c_2$. Utilizing the generalized PL inequality \eqref{eq:polyak_lojasiewicz}, we obtain that \eqref{eq:lyap-increment-2} also holds for all $\tilde{u}\in\mathbb{R}^n$ that satisfy
    \begin{align*}
        \|\tilde{u}-u^\star\|=|\tilde{u}|_{\mathcal{A}_d}\geq \chi(\|e\|):=\mu_1^{-1}\circ\mu_2^{-1}(c_3\|e\|).
    \end{align*}
    Since $\mu_1$ and $\mu_2$ are class $\mathcal{K}_\infty$ functions, the same is true for $\chi$. It follows from \cite[Theorem 3]{sontag2022remarks} that the discere-time system \eqref{eq:discrete-time-gd-with-input} is ISS with respect to $\mathcal{A}_d$. \hfill $\blacksquare$
\end{pf}
\subsection{Input-to-State Stability of Heavy Ball Method}
In this subsection, we prove ISS for the discrete-time system \eqref{eq:discrete-time-system-with-input} when the maps $g$ and $h$ are defined by 
\begin{subequations}\label{eq:discrete-time-heavy-ball}
    \begin{align}
        h((u,w))&=u, \\ g((u,w),f) &= (u-\gamma f + \nu (u-w),u).
    \end{align}
\end{subequations}
where $\gamma\in\mathbb{R}_{+}$ is the \emph{learning rate} and $\nu>0$ is the \emph{momentum parameter}. The choice of \eqref{eq:discrete-time-heavy-ball}  corresponds to \emph{Polyak's Heavy Ball} method \cite{polyak1964some}. Our standing assumption is once again Assumption \ref{asmp:cost_regularity}. We have the following Proposition.

\begin{prop}\label{prop:hb-is-iss}
    Let Assumption \ref{asmp:cost_regularity} be satisfied, and let $\gamma$ and $\nu$ be such that
    \begin{align}\label{eq:stability-condition-hb}
        0&<\gamma<\frac{2}{2+L_\phi}, & 0&<\nu<(4\gamma-2\gamma^2(2+L_\phi))^{\frac{1}{2}}.
    \end{align}
    Then, the discrete time system defined by \eqref{eq:discrete-time-system-with-input} and \eqref{eq:discrete-time-heavy-ball} is ISS with respect to $\mathcal{A}_d=\{(u^\star,u^\star)\}$. 
\end{prop}
\begin{pf}
    We define the Lyapunov function candidate
    \begin{align}\label{eq:lyapunov-function-hb}
        W(\tilde{\eta}):=V(\tilde{u}) + \|\tilde{u}-\tilde{w}\|^2.
    \end{align}
    where $V$ is the Lyapunov function from \eqref{eq:lyapunov-function}. Similar to the proof of Proposition \ref{prop:gd-is-iss}, we obtain from the Descent Lemma \cite[Proposition 6.1.2]{bertsekas2015convex} that
    \begin{align*}
        V(\tilde{u}^+)-V(\tilde{u})\leq \nabla V(\tilde{u})^\top(\tilde{u}^+-\tilde{u})+\frac{L_\phi}{2}\|\tilde{u}^+-\tilde{u}\|^2,
    \end{align*}
    for all $\tilde{u}\in\mathbb{R}^n$. From \eqref{eq:discrete-time-heavy-ball}, we obtain that
    \begin{align*}
        \tilde{u}^+-\tilde{w}^+ = \tilde{u}^+-\tilde{u}.
    \end{align*}
    Hence, with the notation $\Delta W = W(\tilde{\eta}^+)-W(\tilde{\eta})$, we obtain that
    \begin{align*}
        \Delta W\leq \nabla V(\tilde{u})^\top(\tilde{u}^+-\tilde{u})+\tilde{L}_{\phi}\|\tilde{u}^+-\tilde{u}\|^2- \|\tilde{u}-\tilde{w}\|^2,
    \end{align*}
    where $\tilde{L}_{\phi}=1+\frac{L_\phi}{2}$. Substituting from \eqref{eq:discrete-time-heavy-ball} and re-arranging terms, we obtain the inequality
    \begin{align*}
        \Delta W\leq - \kappa(\tilde{\eta})^\top P \kappa(\tilde{\eta}) + c_1\|e\|\|\kappa(\tilde{\eta})\| + c_0 \|e\|^2,
    \end{align*}
    where $\kappa(\tilde{\eta})=(\nabla \phi(\tilde{\eta}),\tilde{u}-\tilde{w})$, $P$ is the block-matrix
    \begin{align*}
        P&= \begin{bmatrix}\gamma (1-\gamma \tilde{L}_{\phi})I & -\frac{\nu}{2}(1-2\gamma \tilde{L}_{\phi})I \\ -\frac{\nu}{2}(1-2\gamma \tilde{L}_{\phi})I & (1-\tilde{L}_{\phi}\nu^2)I \end{bmatrix},
    \end{align*}
    and the constants $c_1,c_0$ are
    \begin{align*}
        c_1&=  \gamma\tilde{L}_{\phi}(4+\nu^2)^{\frac{1}{2}}, & c_0&=\gamma^2\tilde{L}_{\phi}.
    \end{align*}
    From basic Linear Algebra, the matrix $P$ is positive definite if and only if $\det(P)>0$ and $\mathrm{trace}(P)>0$, which leads to the conditions
    \begin{align*}
        \gamma -\gamma^2 \tilde{L}_{\phi} - \frac{\nu^2}{4}&>0, & 1+\gamma - (\gamma^2+\nu^2)\tilde{L}_{\phi}&>0.
    \end{align*}
    both of which holds under the assumption \eqref{eq:stability-condition-hb}. Therefore, when \eqref{eq:stability-condition-hb} is satisfied, the matrix $P$ is positive definite which implies that its smallest eigenvalue is strictly positive, i.e. $\lambda_{\min}(P)>0$. Consequently, we obtain that
    \begin{align*}
        \Delta W\leq -2c_2\|\kappa(\tilde{\eta})\| + c_1\|e\|\|\kappa(\tilde{\eta})\| + c_0 \|e\|^2,
    \end{align*}
    where $c_2=\lambda_{\min}(P)/2$, which implies that
    \begin{align*}
        \Delta W\leq -c_2\|\kappa(\tilde{\eta})\|^2,
    \end{align*}
    for all $\tilde{\eta}\in\mathbb{R}^{2n}$ that satisfy
    \begin{align}\label{eq:lyapunov-increment-hb-2}
        \big(\|\nabla V(\tilde{u})\|+\|\tilde{u}-\tilde{w}\|^2\big)^{\frac{1}{2}}\geq c_3\|e\|,
    \end{align}
    where $c_3$ is once again the positive constant given by $c_3=(c_1^2+4 c_0 c_2)^{\frac{1}{2}}+c_1)/c_2$. Clearly, the function on the left hand side of the inequality \eqref{eq:lyapunov-increment-hb-2} is positive definite with respect to $\mathcal{A}_d=\{(u^\star,u^\star)\}$ and is radially unbounded. It follows from \cite[Theorem 3]{sontag2022remarks} that the discere-time system \eqref{eq:discrete-time-gd-with-input} is ISS with respect to $\mathcal{A}_d$. \hfill $\blacksquare$
\end{pf}
\subsection{Input-to-State Stability of Projected Gradient Descent}
The last example we consider is the discrete-time system
\begin{align}\label{eq:discrete-time-system-pgd-with-input}
    \tilde{u}^+&= \mathcal{P}_{\mathcal{U}}(\tilde{u}-\gamma(\nabla\phi(\tilde{u})+e)),
\end{align}
where $\mathcal{U}$ is a closed and convex set, and $\mathcal{P}_{\mathcal{U}}(\cdot):\mathbb{R}^n\rightarrow \mathcal{U}$ is the projection operator defined by
\begin{align}
    \mathcal{P}_{\mathcal{U}}(u):=\arginf_{\tilde{u}\in\mathcal{U}}\|u-\tilde{u}\|.
\end{align}
The choice \eqref{eq:discrete-time-system-pgd-with-input} corresponds to the Projected Gradient Descent method. To guarantee convergence, we impose the following standard assumption.
\begin{asmp}\label{asmp:strongly-convex-L-smooth}\normalfont
    The set $\mathcal{U}$ is closed and convex, and the function $\phi$ is continuously differentiable, strongly convex with parameter $\mu_\phi>0$, and there exists a constant $L_{\phi}$ such that
    \begin{align}
        \|\nabla\phi(u)-\nabla\phi(\tilde{u})\|\leq L_{\phi}\|u-\tilde{u}\|,
    \end{align}
    for all $u,\tilde{u}\in\mathbb{R}^n$. 
\end{asmp}
Then, we have the following Proposition.

\begin{prop}\label{prop:pgd-is-iss}\normalfont
    Let Assumption \ref{asmp:strongly-convex-L-smooth} be satisfied, and let $\gamma$ be such that $0<\gamma L_\phi<2$. Then, the discrete time system \eqref{eq:discrete-time-system-pgd-with-input} is ISS with respect to $\mathcal{A}_d=\{u^\sharp\}$, where $u^\sharp$ is the unique solution of Problem \eqref{eq:minimization-problem}. 
\end{prop}
\begin{pf}
    First, we observe that, by adding and subtracting a term, the discrete-time system \eqref{eq:discrete-time-system-pgd-with-input} can be re-written as
    \begin{align}
        \tilde{u}^+&= \mathcal{P}_{\mathcal{U}}(\tilde{u}-\gamma\nabla\phi(\tilde{u})) + b(\tilde{u},e),
    \end{align}
    where the vector-valued function $b$ is defined by
    \begin{align*}
        b(\tilde{u},e)&:= \mathcal{P}_{\mathcal{U}}(\tilde{u}-\gamma(\nabla\phi(\tilde{u})+e))-\mathcal{P}_{\mathcal{U}}(\tilde{u}-\gamma\nabla\phi(\tilde{u})).
    \end{align*}
    On the other hand, it is well-known from convex analysis that, when $\mathcal{U}$ is a closed and convex set, the projection $\mathcal{P}_{\mathcal{U}}$ is a non-expansive map, that is 
    \begin{align*}
        \|\mathcal{P}_{\mathcal{U}}(u)-\mathcal{P}_{\mathcal{U}}(\tilde{u})\|\leq \|u-\tilde{u}\|,
    \end{align*}
    for all $u,\tilde{u}\in\mathbb{R}^n$. Hence, for all $\tilde{u},e\in\mathbb{R}^n$, we have that
    \begin{align*}
        \|b(\tilde{u},e)\|\leq \gamma \|e\|.
    \end{align*} 
    To prove ISS, observe that if $\tilde{u}:\mathbb{N}\rightarrow\mathbb{R}^n$ is a solution of the discrete-time system \eqref{eq:discrete-time-system-pgd-with-input} for some input $e:\mathbb{N}\rightarrow\mathbb{R}^n$, then
    \begin{align}
        \tilde{u}(k+1) = \mathcal{T}(u(k))+b(\tilde{u}(k),e(k)),
    \end{align}
    for all $k\in\mathbb{N}$, where $\mathcal{T}(\tilde{u}):= \mathcal{P}_{\mathcal{U}}(\tilde{u}-\gamma\nabla \phi(\tilde{u}))$. Under Assumption \ref{asmp:strongly-convex-L-smooth} and the condition $0<\gamma L_{\phi}<2$, it is a classical result that the map $\mathcal{T}$ is a uniformly contractive map \cite[Proposition 6.1.8]{bertsekas2015convex} with the contraction constant given by $\rho := \max\{|1-\gamma L_{\phi}|,|1-\gamma\mu_{\phi}|\}.$ That is, the map $\mathcal{T}$ has a unique fixed point $u^\sharp$ that coincides with the solution of the minimization problem \eqref{eq:minimization-problem}. 
    We proceed by observing that
    \begin{align*}
        \tilde{u}(k+1)&= \mathcal{T}(\mathcal{T}(\tilde{u}(k-1))+b(\tilde{u}(k-1),e(k-1)))\\ 
        &+ b(\tilde{u}(k),e(k)).
    \end{align*}
    By adding and subtracting a term, we have that
    \begin{align*}
        \tilde{u}(k+1)&= \mathcal{T}\circ\mathcal{T}(\tilde{u}(k-1)) + \tilde{e}_0(k) + \tilde{e}_1(k),
    \end{align*}
    where $\tilde{e}_0$ and $\tilde{e}_1$ are defined by
    \begin{align*}
        \tilde{e}_0(k)&:= b(\tilde{u}(k),e(k)),\\
        \tilde{e}_1(k)&:=\mathcal{T}(\mathcal{T}(\tilde{u}(k-1))+b(\tilde{u}(k-1),e(k-1)))\\
        &- \mathcal{T}\circ\mathcal{T}(\tilde{u}(k-1)).
    \end{align*}
    Since $\mathcal{T}$ is a contraction, it follows that
    \begin{align*}
        \|\tilde{e}_1(k)\|\leq \rho \|b(\tilde{u}(k-1),e(k-1))\|\leq \rho \gamma\|e(k-1)\|.
    \end{align*}
    By recursively applying the steps above, we obtain that
    \begin{align*}
        \tilde{u}(k+1)&= \mathcal{T}^{k+1}(\tilde{u}(0)) + \sum_{i=0}^k \tilde{e}_i(k),
    \end{align*}
    where $\mathcal{T}^k:=\mathcal{T}\circ\mathcal{T}\circ\cdots\mathcal{T}$ is the $k$-times composition of $\mathcal{T}$ with itself, and each $\tilde{e}_i$ satisfies the inequality
    \begin{align*}
        \|\tilde{e}_i(k)\|\leq \rho^i \gamma\|e(k-i)\|.
    \end{align*}
    \begin{figure}[t!]
    \centering
    \includegraphics[width=\linewidth]{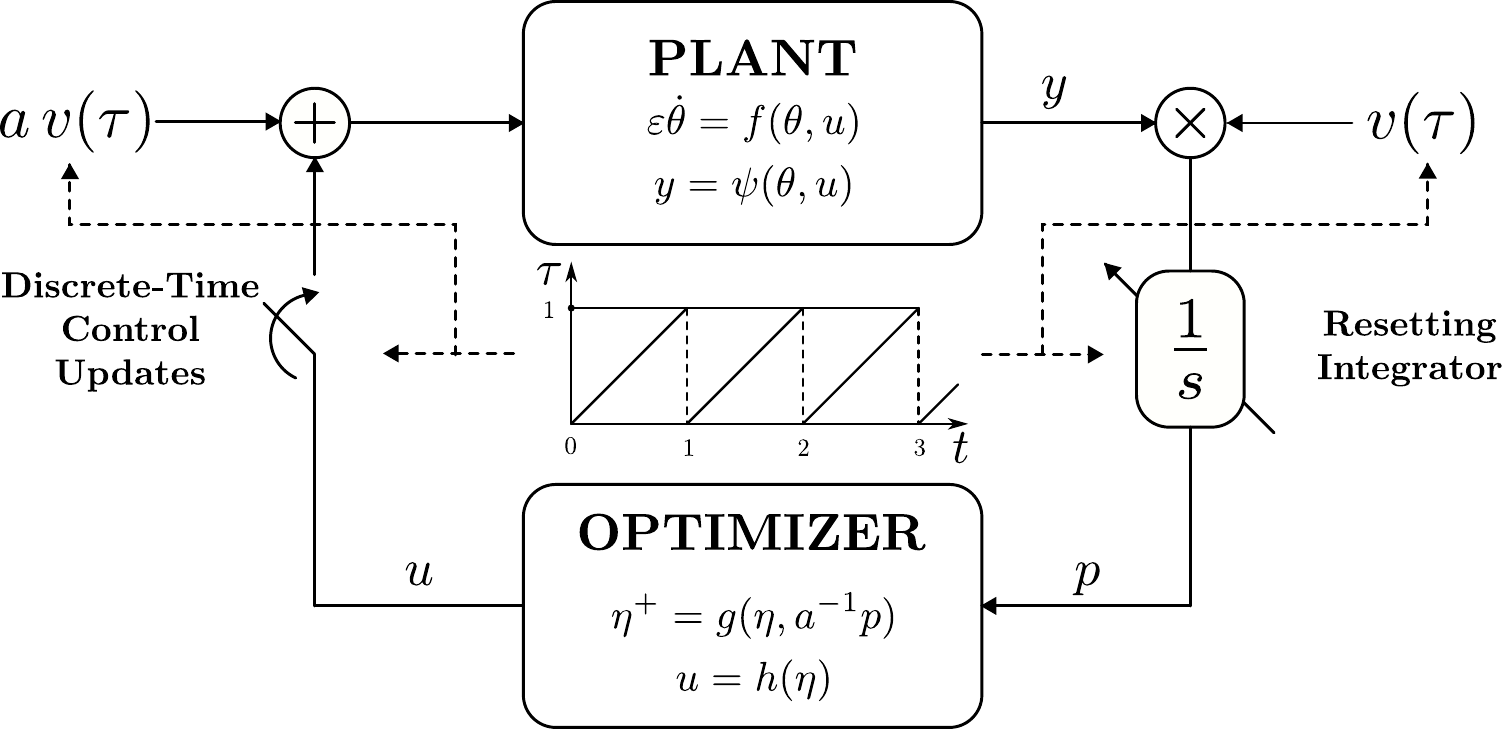}
    \caption{A block diagram description of the interconnection between the proposed model-free optimization algorithm and a dynamic cost function arising as the steady-state output function of a ``fast" plant.}
    \label{fig:block-diagram-dynamic}
\end{figure}
    Through the triangle inequality, we obtain that
    \begin{align*}
        \|\tilde{u}(k+1)-u^\sharp\|&\leq \|\mathcal{T}^{k+1}(\tilde{u}(0))-u^\sharp\|\\
        &+ \gamma \sum_{i=0}^k \rho^i\|e(k-i)\|,
    \end{align*}
    Since, by assumption, $\|e\|_\infty=\sup_{k\in\mathbb{N}}\|e(k)\| < +\infty$, we obtain that
    \begin{align*}
        \|\tilde{u}(k+1)-u^\sharp\|\leq \|\mathcal{T}^{k+1}(\tilde{u}(0))-u^\sharp\| + \gamma\|e\|_\infty \sum_{i=0}^k \rho^i
    \end{align*}
    On the other hand, since $0<\rho<1$, we have that
    \begin{align*}
        \sum_{i=0}^k \rho^i \leq 1/(1-\rho),
    \end{align*}
    for all $k\in\mathbb{N}$, from which we conclude that
    \begin{align*}
        \|\tilde{u}(k)-u^\sharp\|\leq \|\mathcal{T}^{k}(\tilde{u}(0))-u^\sharp\| + \frac{\gamma}{1-\rho}\|e\|_\infty,
    \end{align*}
    for all $k\in\mathbb{N}$. 
    Finally, because $\mathcal{T}(u^\sharp)=u^\sharp$ by definition, and because $\mathcal{T}^k$ is uniformly contractive with a contraction constant $\rho^k$, we obtain that
    \begin{align*}
        \|\tilde{u}(k)-u^\sharp\|\leq \rho^{k}\|\tilde{u}(0)-u^\sharp\| + \frac{\gamma}{1-\rho}\|e\|_\infty,
    \end{align*}
    which is precisely the definition of ISS.  \hfill $\blacksquare$
\end{pf}

\vspace{-0.4cm}
\section{{ Interconnections with Dynamic Plants}} \label{sec:dynamic-plants}
In the previous sections, we established that the proposed hybrid model-free optimization { framework} is capable of inducing UGAS stability of a neighborhood of the optimizer of a \emph{static} cost function under natural assumptions, i.e. Assumption \ref{asmp:cost_regularity} or Assumption \ref{asmp:strongly-convex-L-smooth}. We also provided an explicit characterization of a \emph{bounding box} that contains this UGAS compact set. In the current section, we now show that when the proposed algorithm is interconnected with a \emph{fast} plant, then the same results hold in a semi-global practical sense as the timescale separation between the plant and the optimization algorithm increases.
Consider the continuous-time dynamical system
\begin{align}\label{eq:fast-plant}
    \varepsilon\dot{\theta}&= f(\theta,u), & y&=\psi(\theta,u),
\end{align}
where $\theta\in\mathbb{R}^{l}$ is the state, $u\in\mathbb{R}^n$ is the input, $y\in\mathbb{R}$ is the output, $\psi:\mathbb{R}^{l}\times\mathbb{R}^n\rightarrow\mathbb{R}^n$ is the output map, and $\varepsilon>0$ is a parameter. The following assumption formalizes natural conditions which guarantee that the plant \eqref{eq:fast-plant} is sufficiently well-behaved.
\begin{asmp}\label{asmp:steady-state-output-map}\normalfont
    The following statements holds.
    \begin{enumerate}
        \item The vector field $f$ is Lipschitz continuous;
        \item The map $\psi$ is continuously differentiable;
        \item There exists a continuously differentiable map $\chi:\mathbb{R}^n\rightarrow\mathbb{R}^l$ such that $f(\chi(u),u)=0,$ for all $u\in\mathbb{R}^n$.
        \item $\chi(u)$ is UGAS for \eqref{eq:fast-plant}, for all $u\in\mathbb{R}^n$.
        \item $\chi^{-1}(K)$ is compact for any compact $K\subset\mathbb{R}^n$. 
    \end{enumerate}
\end{asmp}
\begin{figure}
	\centering
	\includegraphics[width=0.75\linewidth]{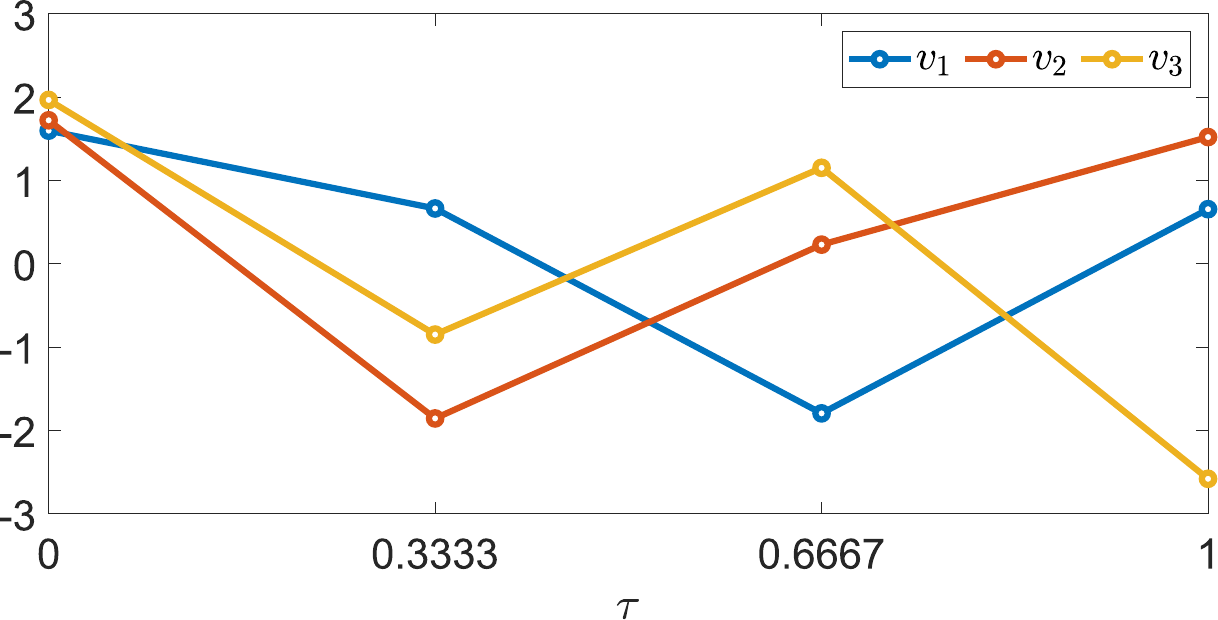}
	\caption{The UPE function defined in \eqref{eq:exploration_signal_example_1}.}
	\label{fig:exploration-signal}
\end{figure}
\begin{figure*}
	\centering
	\includegraphics[width=\linewidth]{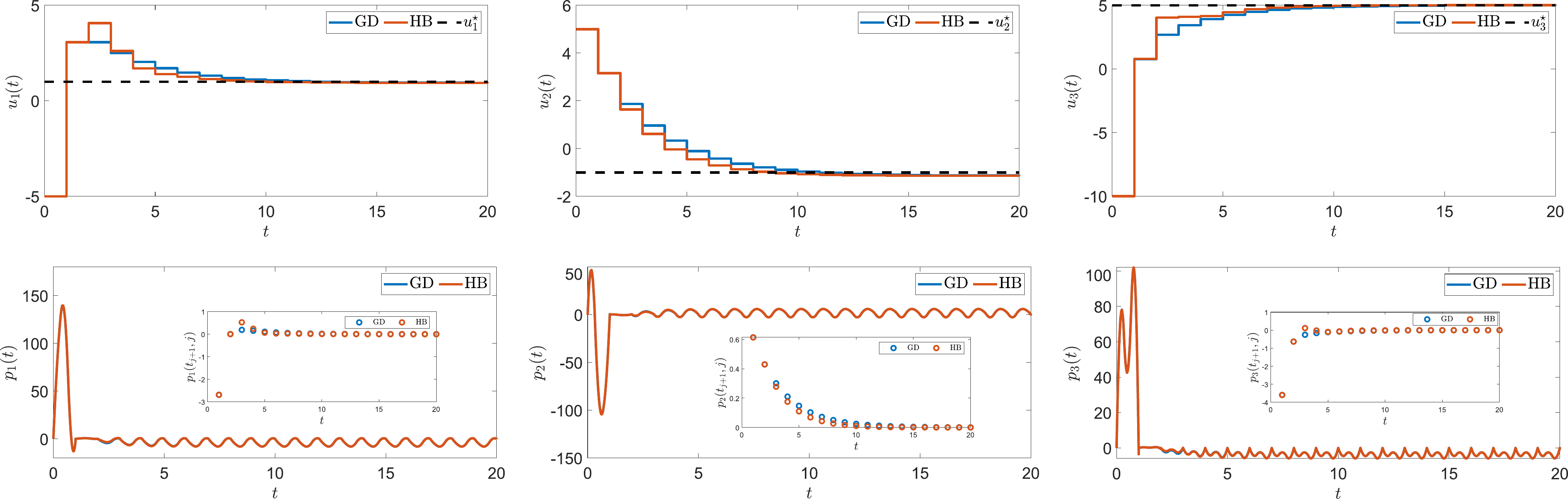}
	\caption{Simulation results for Example 1. The top row shows the evolution of components of the parameter $u$. The bottom row shows the evolution of the components of the state $p$ of the PRLI. The insets in the bottom row show the state $p$ of the PRLI immediately before a jump is triggered, which resets $p$ to $p=0$.}
	\label{fig:example_1_results}
\end{figure*}
Under Assumption \ref{asmp:steady-state-output-map}, the state of the plant \eqref{eq:fast-plant} converges asymptotically to the point $\chi(u)$ if $u$ is a constant, which enables us to define the steady-state output map
\begin{align}\label{eq:steady-state-out-map}
    \phi(u)=\psi(\chi(u),u).
\end{align}
Our goal is to study the behavior of the interconnection between the model-free optimization algorithm represented by the HDS $\mathcal{H}$ and the plant dynamics \eqref{eq:fast-plant} for the purpose of optimizing the function $\phi$. 

To that end, we assume that the function $\phi$ satisfies Assumption \ref{asmp:cost_regularity} or \ref{asmp:strongly-convex-L-smooth}, which would guarantee that, for some compact set $\mathcal{A}\subset\mathbb{R}^{n+m}\times\mathbb{R}^n$, the compact set $\mathcal{A}\times[0,1]$ is UGAS for the HDS $\mathcal{H}$. 

We now investigate the effect of the transient plant dynamics \eqref{eq:fast-plant} on the stability of the HDS $\mathcal{H}$. Specifically, we will study the behavior of the HDS $\hat{\mathcal{H}}_{\varepsilon}=(\mathcal{C},\mathcal{F}_{\varepsilon},\mathcal{D},\mathcal{G})$ where the flow and jump maps are defined by
\begin{align*}
    \mathcal{F}_{\varepsilon}(\eta,\tau,p,\theta)&=(0, v(\tau) y,1,\varepsilon^{-1}f(\theta,h(\eta)+a v(\tau))),\\
    \mathcal{G}(\eta,\tau,p,\theta)&=(g((u,w),a^{-1}p),0,0,\theta),
\end{align*}
$y$ is the output of the plant \eqref{eq:fast-plant}, and the functions $h$ and $g$ are as defined previously. A block diagram depiction of the HDS $\hat{\mathcal{H}}_{\varepsilon}$ is shown in Figure \ref{fig:block-diagram-dynamic}.

Formally, as $\varepsilon\rightarrow 0^+$, a time-scale separation is induced between the fast dynamics of the plant and the optimization algorithm. However, the state-of-the-art results on singular perturbation of HDS \cite{wang2012analysis} necessitate that the fast states are restricted apriori to a compact set. Therefore, to be able to invoke the results in \cite{wang2012analysis}, the flow and jump sets of the interconnected system must be carefully defined. To that end, we let $K\subset \mathbb{R}^{n+m}\times\mathbb{R}^{n}$ be any compact set that contains the set $\mathcal{A}$ in its interior. Next, define the set $K_{\mathcal{U}}\subset\mathbb{R}^n$ by
\begin{align*}
	K_{\mathcal{U}}&= \{u ~|~\exists (w,p) \text{ s.t. }((u,w),p)\in K\}.
\end{align*}
Clearly, the set $K_{\mathcal{U}}$ is compact. In addition, thanks to Assumption \ref{asmp:steady-state-output-map}, we have that the set 
$$\Theta = \chi^{-1}(K_{\mathcal{U}}),$$
is compact. We now define the flow and jump sets for the interconnected HDS $\hat{\mathcal{H}}_{\varepsilon}$ by
\begin{align*}
	\mathcal{C}&= K\times[0,1]\times\Theta, &
	\mathcal{D}&= K\times\{1\}\times\Theta.
\end{align*}
On this form, it is clear that the HDS $\hat{\mathcal{H}}_{\varepsilon}$ belongs to the class of singularly perturbed HDS studied in \cite{wang2012analysis}. Indeed, thanks to Assumption \ref{asmp:steady-state-output-map}, specifically the fact that $\chi(u)$ is UGAS for the plant dynamics \eqref{eq:fast-plant} for every $u\in\mathbb{R}^n$, we conclude that the HDS $\hat{\mathcal{H}}_{\varepsilon}$ has a well-defined \emph{reduced-order} HDS, henceforth denoted by $\mathcal{H}_{rd}$. Moreover, following Example 1 in \cite{wang2012analysis}, we see that the averaged HDS $\mathcal{H}_{rd}=(C_{rd},F,D_{rd},G)$, where
\begin{align*}
	C_{rd}&=K\times[0,1], & D_{rd}&=K\times\{0\},
\end{align*}
and the flow and jump maps $F,G$ coincide with those of the HDS $\mathcal{H}$ as defined in \eqref{eq_SCRLI}, and where $\phi$ is the steady state output map defined in \eqref{eq:steady-state-out-map}.

By construction, we have that $\mathcal{A}\times[0,1]\subset C_{rd}\cup D_{rd}$. Also by construction, we have that $C_{rd}\subset C$ and $D_{rd}\subset D$, where $C$ and $D$ are the flow and jump sets of the HDS $\mathcal{H}$ defined in \eqref{eq_SCRLI}. Since the flow and jump maps of the two HDS coincide, it follows that the HDS $\mathcal{H}_{rd}$ is \emph{contained} in the HDS $\mathcal{H}$ which implies that if the compact set $\mathcal{A}\times[0,1]$ is UGAS for the HDS $\mathcal{H}$ then it is automatically UGAS for the HDS $\mathcal{H}_{rd}$ \cite[Proposition 3.32]{goebel2012hybrid}. However, Theorem \ref{thm:orig-sys-is-ugas} establishes that $\mathcal{A}\times[0,1]$ is indeed UGAS for the HDS $\mathcal{H}$. Thus, the set $\mathcal{A}\times[0,1]$ is UGAS for the averaged HDS $\mathcal{H}_{rd}$. By invoking \cite[Theorem 2]{wang2012analysis}, we obtain that the set  $\mathcal{A}\times[0,1]\times\Theta$ is SGPAS as $\varepsilon\rightarrow 0^+$. The preceeding discussion is the proof of the following Theorem. 
\begin{thm}
    Let Assumption \ref{asmp:steady-state-output-map} be satisfied. Then, there exists a compact set $\mathcal{A}$ such that the set $\mathcal{A}\times[0,1]\times\Theta$ is SGPAS as $\varepsilon\rightarrow 0^+$ for the HDS $\hat{\mathcal{H}}_{\varepsilon}$, provided that:
    \begin{enumerate}
        \item $\phi$ satisfies Assumption \ref{asmp:cost_regularity}, and the functions $g$ and $h$ are given by \eqref{eq:discrete-time-gd} with $0<\gamma L_\phi<2$; OR
        \item $\phi$ satisfies Assumption \ref{asmp:cost_regularity}, and the functions $g$ and $h$ are given by \eqref{eq:discrete-time-heavy-ball} with $\gamma$ and $\nu$ satisfying \eqref{eq:stability-condition-hb}; OR 
        \item $\phi$ satisfies Assumption \ref{asmp:strongly-convex-L-smooth}, and the functions $g$ and $h$ are given by \eqref{eq:discrete-time-system-pgd-with-input} with $e=0$ and $0<\gamma L_\phi<2$.
    \end{enumerate}
\end{thm}
\begin{figure*}
	\centering
	\includegraphics[width=\linewidth]{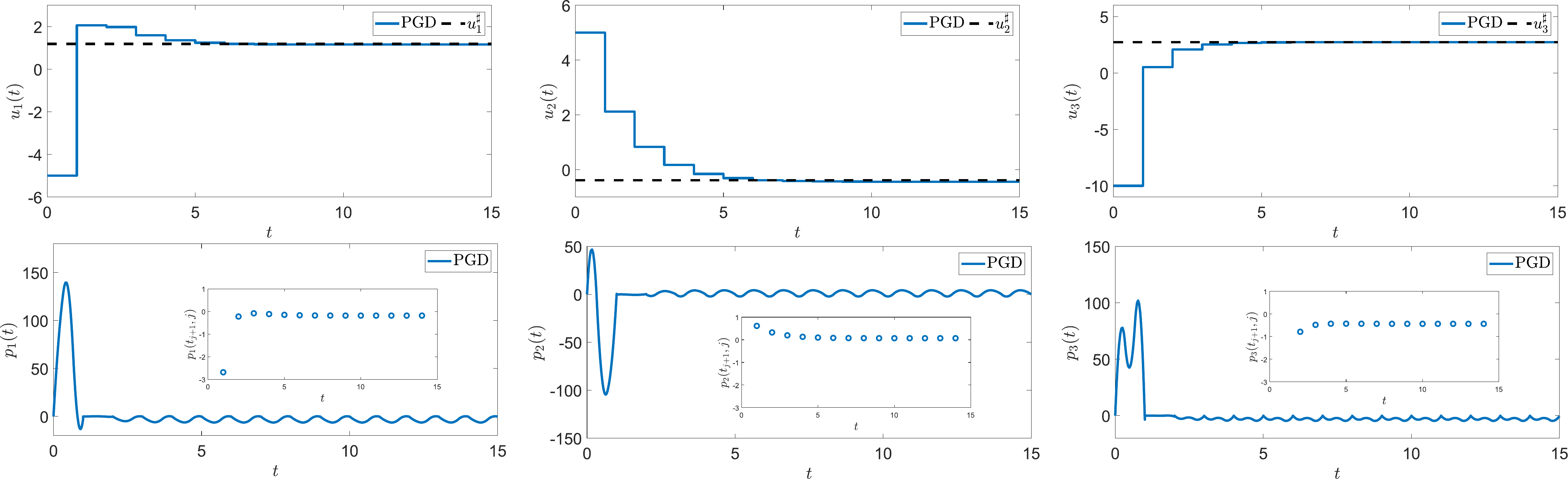}
	\caption{Simulation results for Example 2. The top row shows the evolution of components of the parameter $u$. The bottom row shows the evolution of the components of the state $p$ of the PRLI. The insets in the bottom row show the state $p$ of the PRLI immediately before a jump is triggered, which resets $p$ to $p=0$.}
	\label{fig:example_2_results}
\end{figure*}
\section{Numerical Simulations}\label{sec:numerical-results}
In the current section, we provide several numerical simulations to illustrate the theoretical results presented thus far. Throughout the current section, we take $n=3$ and define the function $v:[0,1]\rightarrow\mathbb{R}^3$ by
\begin{align}\label{eq:exploration_signal_example_1}
    v(\tau):=\begin{cases}
        v^0+3\tau v^1, & \tau\in\left[0,\frac{1}{3}\right),\\
        v^1+(3\tau-1)v^2, & \tau\in\left[\frac{1}{3},\frac{2}{3}\right),\\
        v^2+(3\tau-2)v^3, & \tau\in\left[\frac{2}{3},1\right],
    \end{cases}
\end{align}
where the constants $\{v^0,v^1,v^2,v^3\}\subset\mathbb{R}^3$ are
\begin{align*}
    v^0&\approx (1.5993,1.7234,1.9678), \\
    v^1&\approx (0.6653,-1.8544,-0.8478), \\
    v^2&\approx (-1.7934,0.2314,1.1533), \\
    v^3&\approx (0.6570,1.5225,-2.5788).
\end{align*}
It can be verified that the above function fits Definition \ref{defn:UPE} (up to approximation errors). Moreover, Lemma \ref{lem:UPE-from-PE} can be used to construct a function from $v$ that fits Definition \ref{defn:UPE} exactly, if needed. For illustration purposes, a plot of the function is given in Figure \ref{fig:exploration-signal}.

\begin{figure*}
	\centering
	\includegraphics[width=\linewidth]{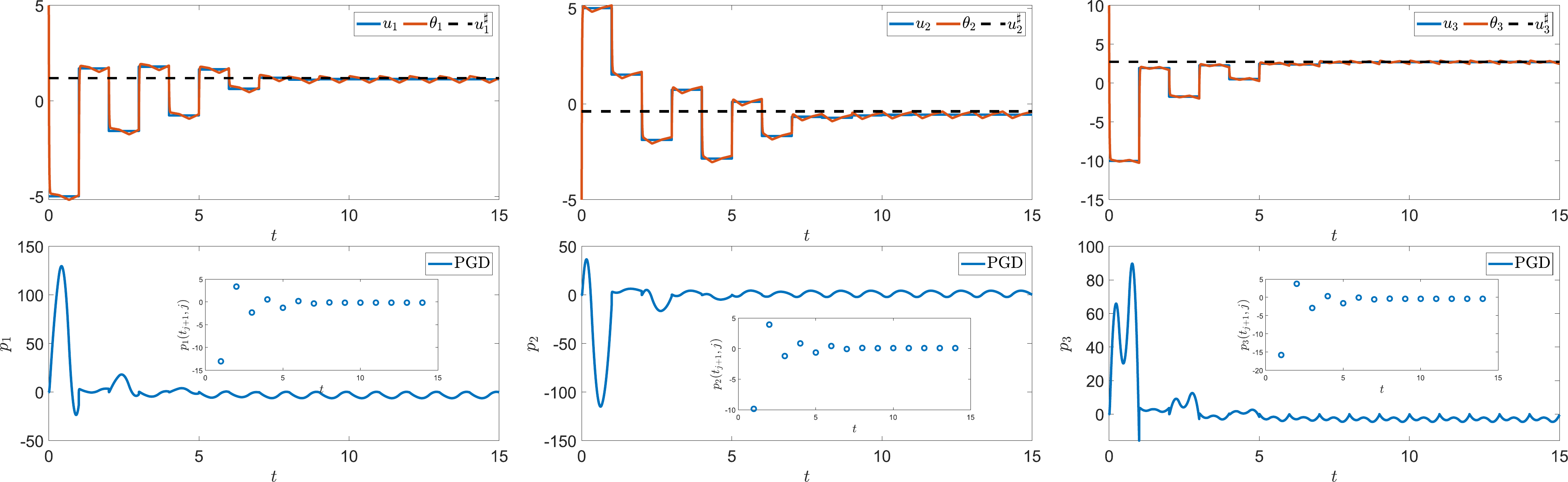}
	\caption{Simulation results for Example 3. The top row shows the evolution of components of the parameter $u$ and the state of the plant $\theta$. The bottom row shows the evolution of the components of the state $p$ of the PRLI. The insets in the bottom row show the state $p$ of the PRLI immediately before a jump is triggered, which resets $p$ to $p=0$.}
	\label{fig:example_3_results}
\end{figure*}

\begin{example}\normalfont
	For the first example, consider the static cost function $\phi:\mathbb{R}^3\rightarrow\mathbb{R}$ defined by
	\begin{subequations}\label{eq:cost_function_example_1}
	\begin{align}
		\phi(u)&=\frac{1}{2} (u-u^\star)^\top Q (u-u^\star) - 20,
	\end{align}
	where the vector $u^\star$ and the symmetric matrix $Q$ are 
	\begin{align}
		u^\star&= \begin{pmatrix} 1 \\ -1 \\ 5\end{pmatrix}, &  Q&=\begin{pmatrix} 2 & 0 & 1\\ 0 & 1 & 0\\ 1 & 0 & 2\end{pmatrix}.
	\end{align}
	\end{subequations}
	Clearly, the function $\phi$ satisfies Assumption \ref{asmp:cost_regularity}. We simulate the HDS $\mathcal{H}$ for the two cases  when the functions $g$ and $h$ correspond to gradient descent, i.e. \eqref{eq:discrete-time-gd}, and the heavy ball method, i.e. \eqref{eq:discrete-time-heavy-ball}. In both cases, we take $a=0.1$ and $\gamma = 0.5$, and, for the heavy ball method, we use $\nu = 0.125$ as the value of the momentum parameter. Simulation results are shown in Figure \ref{fig:example_1_results}. It is readily observed that both algorithms converge to the minimizer of the function $\phi$ and that the convergence of the heavy ball method is faster than gradient descent as expected.
\end{example}

\begin{example}\normalfont
	Consider once again the static cost function $\phi:\mathbb{R}^3\rightarrow\mathbb{R}$ defined by \eqref{eq:cost_function_example_1} and let the function $v:[0,1]\rightarrow\mathbb{R}^3$ be defined by \eqref{eq:exploration_signal_example_1}. Define the set
    \begin{align}\label{eq:example_2_feasible_set}
        \mathcal{U}&= \{u\in\mathbb{R}^3~|~\|u\|^2\leq 9\}.
    \end{align}
    Clearly, $\phi$ and $\mathcal{U}$ satisfy Assumption \ref{asmp:strongly-convex-L-smooth}. Therefore, $\phi$ has a unique minimizer $u^\sharp$ over the feasible set $\mathcal{U}$. By applying standard techniques from convex optimization \cite{bertsekas2015convex}, it can be shown that the unique minimizer $u^\sharp$ of $\phi$ over the feasible set $\mathcal{U}$ is approximately given by 
	\begin{align}
		u^\sharp \approx (1.18,-0.38, -2.72).
	\end{align}
	We now simulate the HDS $\mathcal{H}$ for the case when the functions $g$ and $h$ correspond to the projected gradient descent. We take $a=0.1$ and $\gamma = 0.5$. Simulation results are shown in Figure \ref{fig:example_2_results}. It is readily observed that the projected gradient descent method converges to the unique minimizer $u^\sharp$. Moreover, by construction of the flow and jump maps, the parameter $u$ remains within the feasible set $\mathcal{U}$ for all time (beyond the first jump).
\end{example}

\begin{example}\normalfont
    In the last example, we consider the situation wherein the cost function $\phi$ arises as the steady-state output map of a dynamical system. For simplicity, we consider the dynamical system
    \begin{align}\label{eq:example_2_plant_dynamics}
        \varepsilon \dot{\theta}&=  u - \theta, & y&= \frac{1}{2}(\theta-u^\star)^\top Q (\theta-u^\star)-20,
    \end{align}
    where $\varepsilon > 0$ and, as before, $u^\star$ and $Q$ are as defined in \eqref{eq:cost_function_example_1}. We also define the set $\mathcal{U}\subset\mathbb{R}^3$ by \eqref{eq:example_2_feasible_set} and we let the function $v:[0,1]\rightarrow\mathbb{R}^3$ be defined by \eqref{eq:exploration_signal_example_1}. It is readily seen that the quasi-steady state map $\chi$ is given by $\chi(u)=u$, and the plant \eqref{eq:example_2_plant_dynamics} satisfies Assumption \ref{asmp:steady-state-output-map}. From the definition of the steady-state-output-map \eqref{eq:steady-state-out-map}, we see that $\phi$ is as defined in \eqref{eq:cost_function_example_1}.  Similar to Example 2, the function $\phi$ and the set $\mathcal{U}$ satisfy Assumption \ref{asmp:strongly-convex-L-smooth}, and so the function $\phi$ has the same unique minimizer $u^\sharp$ over the feasible set $\mathcal{U}$. We now simulate the singularly perturbed HDS $\hat{\mathcal{H}}_{\varepsilon}$ for the case when the functions $g$ and $h$ correspond to the projected gradient descent method. We take $a=0.1$, $\gamma = 0.5$, and $\varepsilon = 0.01$. Simulation results are shown in Figure \ref{fig:example_3_results}.
\end{example}


\vspace{-0.3cm}
\section{Conclusion}\label{sec:conclusion}

\vspace{-0.2cm}
We introduced a novel framework for model-free zeroth-order feedback optimization  based on PRLIs. The proposed approach is modular and aims to synergistically integrate the literature on continuous-time optimization of dynamic plants with the rich literature on discrete-time optimization algorithms.Under natural assumptions on the cost function, global asymptotic stability of a compact neighborhood around the minimizer is established for the static cost case. { In contrast to the majority of model-free real-time zeroth-order optimization schemes}, an explicit characterization of this UGAS set is provided without any tuning of ``small" parameters. The results are then extended to the case of a dynamic cost using timescale separation. Our results open the door to several potential future contributions, including distributed frameworks of gradient estimation in networked systems using PRLIs with \emph{asynchronous} resets, as well as the incorporation of stochastic learning and optimization mechanisms into the feedback loops.

\bibliographystyle{elsarticle-num}  
\bibliography{Bibliography}

\begin{thebibliography}{10}
\expandafter\ifx\csname url\endcsname\relax
  \def\url#1{\texttt{#1}}\fi
\expandafter\ifx\csname urlprefix\endcsname\relax\def\urlprefix{URL }\fi
\expandafter\ifx\csname href\endcsname\relax
  \def\href#1#2{#2} \def\path#1{#1}\fi

\bibitem{ghaffari2014power}
A.~Ghaffari, M.~Krsti{\'c}, S.~Seshagiri, Power optimization and control in
  wind energy conversion systems using extremum seeking, IEEE transactions on
  control systems technology 22~(5) (2014) 1684--1695.

\bibitem{colombino2019online}
M.~Colombino, E.~Dall’Anese, A.~Bernstein, Online optimization as a feedback
  controller: Stability and tracking, IEEE Transactions on Control of Network
  Systems 7~(1) (2019) 422--432.

\bibitem{yu2021extremum}
H.~Yu, S.~Koga, T.~R. Oliveira, M.~Krstic, Extremum seeking for traffic
  congestion control with a downstream bottleneck, Journal of Dynamic Systems,
  Measurement, and Control 143~(3) (2021) 031007.

\bibitem{guay2004adaptive}
M.~Guay, D.~Dochain, M.~Perrier, Adaptive extremum seeking control of
  continuous stirred tank bioreactors with unknown growth kinetics, Automatica
  40~(5) (2004) 881--888.

\bibitem{wang1999optimizing}
H.-H. Wang, M.~Krsti{\'c}, G.~Bastin, Optimizing bioreactors by extremum
  seeking, International Journal of Adaptive Control and Signal Processing
  13~(8) (1999) 651--669.

\bibitem{abdelgalil2023singularly}
M.~Abdelgalil, A.~Eldesoukey, H.~Taha, Singularly perturbed averaging with
  application to bio-inspired 3d source seeking, in: 2023 American Control
  Conference (ACC), IEEE, 2023, pp. 885--890.

\bibitem{hauswirth2024optimization}
A.~Hauswirth, Z.~He, S.~Bolognani, G.~Hug, F.~D{\"o}rfler, Optimization
  algorithms as robust feedback controllers, Annual Reviews in Control 57
  (2024) 100941.

\bibitem{scheinker2024100}
A.~Scheinker, 100 years of extremum seeking: A survey, Automatica 161 (2024)
  111481.

\bibitem{Leblanc1922}
M.~Leblanc, Sur l’electri”cation des chemins de fer au moyen de courants
  alternatifs de frequence elevee, Revue Generale de l’Electricite (1922).

\bibitem{krstic2000stability}
M.~Krsti{\'c}, H.-H. Wang, Stability of extremum seeking feedback for general
  nonlinear dynamic systems, Automatica 36~(4) (2000) 595--601.

\bibitem{ariyur2003real}
K.~B. Ariyur, M.~Krstic, Real-time optimization by extremum-seeking control,
  John Wiley \& Sons, 2003.

\bibitem{tan2010extremum}
Y.~Tan, W.~H. Moase, C.~Manzie, D.~Ne{\v{s}}i{\'c}, I.~M. Mareels, Extremum
  seeking from 1922 to 2010, in: Proceedings of the 29th Chinese control
  conference, IEEE, 2010, pp. 14--26.

\bibitem{oliveira2016extremum}
T.~R. Oliveira, M.~Krsti{\'c}, D.~Tsubakino, Extremum seeking for static maps
  with delays, IEEE Transactions on Automatic Control 62~(4) (2016) 1911--1926.

\bibitem{zhu2022extremum}
Y.~Zhu, E.~Fridman, Extremum seeking via a time-delay approach to averaging,
  Automatica 135 (2022) 109965.

\bibitem{yilmaz2024prescribed}
C.~T. Yilmaz, M.~Krstic, Prescribed-time extremum seeking for delays and pdes
  using chirpy probing, IEEE Transactions on Automatic Control (2024).

\bibitem{khong2013unified}
S.~Z. Khong, D.~Ne{\v{s}}i{\'c}, Y.~Tan, C.~Manzie, Unified frameworks for
  sampled-data extremum seeking control: Global optimisation and multi-unit
  systems, Automatica 49~(9) (2013) 2720--2733.

\bibitem{hazeleger2022sampled}
L.~Hazeleger, D.~Ne{\v{s}}i{\'c}, N.~van~de Wouw, Sampled-data extremum-seeking
  framework for constrained optimization of nonlinear dynamical systems,
  Automatica 142 (2022) 110415.

\bibitem{khong2015extremum}
S.~Z. Khong, Y.~Tan, C.~Manzie, D.~Ne{\v{s}}i{\'c}, Extremum seeking of
  dynamical systems via gradient descent and stochastic approximation methods,
  Automatica 56 (2015) 44--52.

\bibitem{durr2013lie}
H.-B. D{\"u}rr, M.~S. Stankovi{\'c}, C.~Ebenbauer, K.~H. Johansson, Lie bracket
  approximation of extremum seeking systems, Automatica 49~(6) (2013)
  1538--1552.

\bibitem{abdelgalil2022recursive}
M.~Abdelgalil, H.~Taha, Recursive averaging with application to bio-inspired
  3-d source seeking, IEEE Control Systems Letters 6 (2022) 2816--2821.

\bibitem{poveda2017framework}
J.~I. Poveda, A.~R. Teel, A framework for a class of hybrid extremum seeking
  controllers with dynamic inclusions, Automatica 76 (2017) 113--126.

\bibitem{abdelgalil2025lie}
M.~Abdelgalil, J.~I. Poveda, On {L}ie-bracket averaging for hybrid dynamical
  systems with applications to model-free control and optimization, IEEE
  Transactions on Automatic Control (2025).

\bibitem{abdelgalil2023multi}
M.~Abdelgalil, D.~E. Ochoa, J.~I. Poveda, Multi-time scale control and
  optimization via averaging and singular perturbation theory: From {ODE}s to
  hybrid dynamical systems, Annual Reviews in Control 56 (2023) 100926.

\bibitem{teel1999semi}
A.~R. Teel, J.~Peuteman, D.~Aeyels, Semi-global practical asymptotic stability
  and averaging, Systems \& control letters 37~(5) (1999) 329--334.

\bibitem{abdelgalil2024initialization}
M.~Abdelgalil, J.~I. Poveda, Initialization-free {L}ie-bracket extremum
  seeking, Systems \& Control Letters 191 (2024) 105881.

\bibitem{goebel2012hybrid}
R.~Goebel, R.~G. Sanfelice, A.~R. Teel, Hybrid Dynamical Systems: Modeling,
  Stability, and Robustness, Princeton University Press, 2012.

\bibitem{subbaraman2016equivalence}
A.~Subbaraman, A.~R. Teel, On the equivalence between global recurrence and the
  existence of a smooth lyapunov function for hybrid systems, Systems \&
  Control Letters 88 (2016) 54--61.

\bibitem{sastry2011adaptive}
S.~Sastry, M.~Bodson, Adaptive control: stability, convergence and robustness,
  Courier Corporation, 2011.

\bibitem{polyak1987introduction}
B.~T. Polyak, Introduction to Optimization, Optimization Software Inc., 1987.

\bibitem{nestruev2003smooth}
J.~Nestruev, Smooth {M}anifolds and {O}bservables, Vol. 220, Springer, 2020.

\bibitem{poveda2017robust}
J.~I. Poveda, A.~R. Teel, A robust event-triggered approach for fast
  sampled-data extremization and learning, IEEE Transactions on Automatic
  Control 62~(10) (2017) 4949--4964.

\bibitem{lessard2022analysis}
L.~Lessard, The analysis of optimization algorithms: A dissipativity approach,
  IEEE Control Systems Magazine 42~(3) (2022) 58--72.

\bibitem{khalil2002nonlinear}
H.~Khalil, Nonlinear Systems, Prentice Hall, 2002.

\bibitem{sontag2022remarks}
E.~D. Sontag, Remarks on input to state stability of perturbed gradient flows,
  motivated by model-free feedback control learning, Systems \& Control Letters
  161 (2022) 105138.

\bibitem{bertsekas2015convex}
D.~Bertsekas, Convex optimization algorithms, Athena Scientific, 2015.

\bibitem{polyak1964some}
B.~T. Polyak, Some methods of speeding up the convergence of iteration methods,
  Ussr computational mathematics and mathematical physics 4~(5) (1964) 1--17.

\bibitem{wang2012analysis}
W.~Wang, A.~R. Teel, D.~Ne{\v{s}}i{\'c}, Analysis for a class of singularly
  perturbed hybrid systems via averaging, Automatica 48~(6) (2012) 1057--1068.

\end{thebibliography}

\end{document}